\documentclass{article}
\usepackage{float, graphicx}
\usepackage[]{epsfig}
\usepackage{amsmath, amsthm, amssymb,}
\usepackage{epsfig}
\usepackage{verbatim}
\usepackage{multicol}
\usepackage{url}
\usepackage[ruled,vlined]{algorithm2e}
\usepackage{latexsym}
\usepackage{mathrsfs}
\usepackage[colorlinks, bookmarks=true]{hyperref}
\usepackage{graphicx}
\usepackage{amsmath}
\usepackage{enumerate}
\usepackage[normalem]{ulem}
\usepackage{bm}
\usepackage{dsfont}
\usepackage{stmaryrd}
\usepackage{enumitem}
\usepackage{tikz}
\usepackage{mathtools}
\usepackage{authblk}
\usepackage[noadjust]{cite}
\usepackage{color}
\usepackage{xcolor}

\usepackage{authblk}
\usepackage{geometry}
\numberwithin{equation}{section}

\renewcommand{\cal}{\mathcal}
\newcommand\cA{{\mathcal A}}
\newcommand\cB{{\mathcal B}}
\newcommand{\cC}{{\cal C}}

\newcommand{\fc}{{\mathfrak c}}

\newcommand{\bma}{{\bm{a}}}
\newcommand{\bmb}{{\bm{b}}}

\newcommand{\bmu}{{\bm{u}}}
\newcommand{\bmv}{{\bm{v}}}
\newcommand{\bmw}{{\bm{w}}}

\newcommand{\bmS}{\bm S}

\newcommand{\bmX}{\bm X}
\newcommand{\bmW}{\bm W}

\newcommand{\rd}{{\rm d}}

\newcommand{\ri}{\mathrm{i}}

\newcommand{\bC}{{\mathbb C}}

\newcommand{\bE}{\mathbb{E}}

\newcommand{\bP}{\mathbb{P}}

\newcommand{\bR}{{\mathbb R}}

\newcommand{\bI}{\mathbb{I}}

\newcommand{\la}{\lambda}

\DeclareMathOperator{\Tr}{Tr}

\DeclareMathOperator{\OO}{O}
\DeclareMathOperator{\oo}{o}

\renewcommand{\Re}{\mathop{\mathrm{Re}}}
\renewcommand{\Im}{\mathop{\mathrm{Im}}}

\newcommand{\deq}{\mathrel{\mathop:}=} 
 
\renewcommand{\leq}{\leqslant}
\renewcommand{\geq}{\geqslant}
\newcommand{\floor}[1] {\lfloor {#1} \rfloor}
\newcommand{\ceil}[1]  {\lceil  {#1} \rceil}

\newcommand{\del}{\partial}

\newcommand{\qq}[1]{[\![{#1}]\!]}

\newcommand{\beq}{\begin{equation}}
\newcommand{\eeq}{\end{equation}}

\newcommand{\Mat}{\mathsf{Mat}}

\theoremstyle{plain} 
\newtheorem{theorem}{Theorem}[section]
\newtheorem*{theorem*}{Theorem}
\newtheorem{lemma}[theorem]{Lemma}
\newtheorem*{lemma*}{Lemma}

\newtheorem*{corollary*}{Corollary}

\newtheorem*{proposition*}{Proposition}
\newtheorem{assumption}[theorem]{Assumption}
\newtheorem*{assumption*}{Assumption}

\newtheorem*{definition*}{Definition}

\newtheorem*{example*}{Example}
\newtheorem{remark}[theorem]{Remark}

\newtheorem*{remark*}{Remark}
\newtheorem*{remarks*}{Remarks}

  \title{Long Random Matrices and Tensor Unfolding}
  
 \author[1]{G{\' e}rard~Ben Arous\thanks{benarous@cims.nyu.edu}}
    \author[2]{Daniel Zhengyu~Huang\thanks{dzhuang@caltech.edu}}

  \author[3]{Jiaoyang~Huang\thanks{jh4427@nyu.edu}}

\affil[1,3]{Courant Institute, NYU, New York, NY}
\affil[2]{California Institute of Technology, Pasadena, CA}

\date{}

\begin{document}



\maketitle

\begin{abstract}
In this paper, we consider the singular values and singular vectors of low rank perturbations of large rectangular random matrices, in the regime the matrix is ``long":
we allow the number of rows (columns) to grow polynomially in the number of columns (rows). We prove there exists a critical signal-to-noise ratio (depending on the dimensions of the matrix), and the extreme singular values and singular vectors exhibit a BBP type phase transition. As a main application, we investigate the tensor unfolding algorithm for the asymmetric rank-one spiked tensor model, and obtain an exact threshold, which is independent of the procedure of tensor unfolding. If the signal-to-noise ratio is above the threshold, tensor unfolding detects the signals; otherwise, it fails to capture the signals.   
\end{abstract}

\section{Introduction}

High order arrays, or tensors have been actively considered in neuroimaging analysis, topic modeling, signal processing and recommendation system \cite{frolov2017tensor, Comon,Hack,Kar,Rendle,zhou, simony2016dynamic,cichocki2015tensor,sidiropoulos2017tensor}.
 Researchers have made tremendous efforts to innovate effective methods for the analysis of tensor data. 
 The \emph{spiked tensor model}, introduced in \cite{NIPS2014_5616} by Richard and Montanari, captures a number of statistical estimation tasks in which we need to extract information from a noisy high-dimensional data tensor. We are given a tensor $\bmX\in (\bR^n)^{\otimes k}$ in the following form
\begin{align*}
\bmX=\beta \bmv_1\otimes \bmv_2\otimes \cdots \otimes \bmv_k+\bmW,
\end{align*}
where $\bmW$ is a noise tensor, $\beta>0$ corresponds to the signal-to-noise ratio (SNR), and 
$\bmv_1\otimes \bmv_2\otimes \cdots \otimes \bmv_k$ is a rank one unseen signal tensor to be recovered.

When $k=2$, the spiked tensor model reduces to the \emph{spiked matrix model} of the form ``signal $+$ noise", which has been intensively studied in the past twenty years. 
It is now well-understood that the extreme eigenvalues of low rank perturbations of large random matrices undergo a so-called BBP phase transition \cite{baik2005phase} along
with the change of signal-to-noise ratio, first discovered by Baik, P{\' e}ch{\' e} and the first author. There is an order one critical signal-to-noise ratio $\beta_{\rm c}$, such that below $\beta_{\rm c}$, it is information-theoretical impossible to detect the spikes
\cite{perry2016optimality,onatski2013asymptotic,montanari2015limitation}, and above $\beta_{\rm c}$, it is possible to detect the
spikes by Principal Component Analysis (PCA). A body of work has quantified the behavior of PCA in this setting \cite{johnstone2001distribution, baik2005phase,baik2006eigenvalues,paul2007asymptotics, benaych2012singular,bai2012sample, johnstone2009consistency, birnbaum2013minimax,cai2013sparse,ma2013sparse,vu2013minimax,cai2015optimal,el2008spectrum,ledoit2012nonlinear,donoho2018optimal}.
We refer readers to the review article \cite{johnstone2018pca} by Johnstone and Paul,
for more discussion and references to this and related lines of work.

For the spiked tensor model with $k\geq 3$, the same as the spiked matrix model,  
there is an order one critical signal-to-noise ratio $\beta_{k}$ (depending on the order $k$), such that below $\beta_{k}$, it is information-theoretical impossible to detect the spikes, and above $\beta_{k}$, the maximum likelihood estimator is a distinguishing statistics \cite{chen2019phase,chen2018phase,lesieur2017statistical,perry2020statistical,jagannath2018statistical}.
In the matrix setting the
maximum likelihood estimator is the top eigenvector, which can be computed in polynomial time
by, e.g., power iteration. However, for order $k\geq 3$ tensor, computing the maximum likelihood estimator is NP-hard in generic setting. In this setting, two ``phase transitions" are often studied. There is the critical signal-to-noise ratio SNR$_{stat}=\beta_k$, below which  it is statistically impossible to estimate the parameters. Although above the threshold SNR$_{stat}$, it is possible to estimate the parameters in theory,  there is no known efficient algorithm (polynomial time) to achieve recovery close to the statistical threshold SNR$_{stat}$. Thus many algorithm
development fields  are interested in another critical threshold SNR$_{comp}\geq {}$SNR$_{stat}$ below which it is impossible for an efficient algorithm to achieve recovery. For the spiked tensor model with $k\geq 3$, it is widely believed there exists a computational-to-statistical gaps SNR$_{comp}- {}$SNR$_{stat}>0$.  We refer readers to the article \cite{bandeira2018notes} by Bandeira, Perry and Wein, for more detailed discussion on this phenomenon.


In the work \cite{NIPS2014_5616} of Richard and Montanari, the algorithmic aspects of the spiked tensor model have been studied.  They showed that  tensor power iteration and approximate message passing algorithm with random initialization recovers the signal provided $\beta\gtrsim n^{(k-1)/2}$.  Based on heuristic arguments, they predicted that the necessary and sufficient 
condition for power iteration and approximate message passing (AMP)  algorithm to succeed is $\beta\gtrsim n^{(k-2)/2}$.
This threshold was proven in \cite{lesieur2017statistical} for AMP by Lesieur, Miolane, Lelarge, Krzakala and Zdeborov{\'a}, for power iteration  by Yang, Cheng and the last two authurs \cite{huang2020power}. The same threshold was also achieved by using  gradient descent and Langevin dynamics as studied by Gheissari, Jagannath and the first author \cite{arous2020algorithmic}.

In \cite{NIPS2014_5616}, Richard and Montanari also proposed a method based on tensor unfolding, which unfolds the tensor $\bmX$ to an $n^q\times n^{k-q}$ matrix $\Mat(\bmX)$ for some $1\leq q\leq k-1$,
\begin{align}\label{e:MatX}
\Mat(\bmX)=\beta \bmv\bmu^\top+Z\in \bR^{n^q}\times \bR^{n^{k-q}}.
\end{align}
By taking $q=\floor{k/2}$ and performing matrix PCA on $\Mat(\bmX)$, they proved that the tensor unfolding algorithm reovers the signal provided $\beta\gtrsim n^{(\lceil k/2\rceil-1)/2}$, and predicted that the necessary and sufficient condition  for tensor unfolding  is $\beta\gtrsim n^{(k-2)/4}$.

Several other sophisticated algorithms for the spiked tensor model have been investigated in literature which achieves the sharp threshold $\beta\gtrsim n^{(k-2)/4}$: Sum-of-Squares algorithms \cite{hopkins2015tensor, hopkins2016fast, kim2017community},  sophisticated iteration algorithms \cite{luo2020sharp, zhang2018tensor,han2020optimal}, and an averaged version of gradient descent \cite{biroli2020iron} by Biroli, Cammarota, Ricci-Tersenghi.
The necessary part of this threshold still remains open. Its relation with hypergraphic planted clique problem was discussed in \cite{luo2020open,luo2020tensor} by Luo and Zhang. Its proven for $k=3$ in \cite{hopkins2015tensor} by Hopkins, Shi and Steurer, degree-$4$ sum-of-squares relaxations fail below this threshold.
The candscape complexity of spiked tensor model was studied in \cite{arous2019landscape} by Mei, Montanari, Nica and the first author. A new framework based on the Kac-Rice method that
allows  to compute the annealed complexity of the landscape has been proposed in \cite{ros2019complex} by Ros, Biroli,Cammarota and the first author, which was later used to analyze gradient-based algorithms in non-convex setting \cite{sarao2019afraid,mannelli2020complex} by Mannelli, Biroli,  Cammarota, Krzakala, Urbani, and Zdeborov{\'a}.

In this paper, we revisit the tensor unfolding algorithm introduced by Richard and Montanari. The unfolded matrix $\Mat(\bmX)$ from \eqref{e:MatX} is a spiked matrix model in the form of ``signal $+$ noise". However, it is different from spiked matrix models in random matrix literature, which requires the dimensions of the matrix to be comparable, namely the ratio of the number of rows and the number of columns converges to a fixed constant. In this setting, the singular values and singular vectors of the spiked matrix model \eqref{e:MatX} has been studied in \cite{ benaych2012singular} by Benaych-Georges and Nadakuditi. 

For the unfolded matrix $\Mat(\bmX)$, its dimensions are not comparable. As the size $n$ of the tensor goes to infinity, the ratio of the number of rows and columns goes to zero or infinity, unless $q=k/2$ (in this case $\Mat(\bmX)$ is a square matrix). In this paper, we study the singular values and singular vectors of the spiked matrix model \eqref{e:MatX} in the case where the number of rows (columns) grows polynomially in the number of columns (rows), which we call \emph{low rank perturbation of long random matrices}. In the case when $\beta=0$, the estimates of singular values and singular vectors for long random matrices follow from \cite{alex2014isotropic} by Alex, Erd{\H o}s, Knowles, Yau, and Yin.

To study the low rank perturbations of long random matrices, we use the master equations from \cite{ benaych2012singular}, which characterize the outliers of the perturbed random matrices, and the associated singular values. To analyze the master equations, we use the estimates of singular values and Green's functions of long random matrices from \cite{alex2014isotropic}. Comparing with the setting that the ratio of the number of rows and the number of columns converges to a fixed constant, the challenge is to obtain uniform estimates for the errors in the master equations, which depends only on the number of rows. In this way, we can allow the number of columns to grow much faster than the number of rows. 
For the low rank perturbation of long random matrices, we prove there exists a critical signal-to-noise ratio $\beta_{\rm c}$ (depending on the dimensions of the matrix), and it exhibit a BBP type phase transition. 
We also obtain estimates of the singular values and singular vectors for this model. 
Moreover, we also have precise estimates when the signal-to-noise ratio $\beta$ is close to the threshold $\beta_{\rm c}$. In particular, our results also apply when $\beta$ is close to $\beta_{\rm c}$ in mesoscopic scales. In an independent work \cite{feldman2021spiked} by Feldman, this model has been studied under different assumptions. We refer to Section \ref{s:lowrank} for a more detailed discussion of the differences.

Our results for low rank perturbation of long random matrices can be used to study the unfolded matrix $\Mat(\bmX)$. For the signal-to-noise ratio $\beta=\la n^{(k-2)/4}$ and any $1\leq q\leq k-1$, if $\la>1$, the PCA on $\Mat(\bmX)$ detects the signal tensor; if $\la<1$, the PCA on $\Mat(\bmX)$ fails to capture the signal tensor. This matches the conjectured algorithmic threshold for the spiked tensor model from \cite{NIPS2014_5616}. It is worth mentioning that the threshold we get is independent of the tensor unfolding procedure, namely, it is independent of the choice of $q$. For $q>1$, a further recursive unfolding is needed to recover individual signals $\bmv_i$, which increases the computational cost. We propose to simply take $q=1$ in the tensor unfolding algorithm for each coordinate axis, and unfold the tensor to an $n\times n^{k-1}$ matrix, which gives good approximation of  individual signals $\bmv_i$, provided $\la>1$. In tensor literature, this algorithm is exactly the truncated higher order singular value decomposition (HOSVD) introduced in \cite{de2000multilinear} by De Lathauwer, De Moor and Vandewalle. Later, they developed the higher order orthogonal iteration (HOOI) in \cite{de2000best}, which uses the truncated HOSVD as initialization combining with a power iteration, to find the best low-multilinear-rank approximation of a tensor. 
The performance of HOOI was analyzed in \cite{zhang2018tensor} by Zhang and Xia, for the spiked tensor model. It was proven that for the signal-to-noise ratio $\beta=\la n^{(k-2)/4}$ satisfies $\lambda\geq C_{gap}$ for some large constant $C_{gap}>0$, HOOI converges within a logarithm factor of iterations.

The paper is organized as follows. In Section \ref{s:main}, we state the main results on the singular values and vectors of low rank perturbations of long random matrices.
In Section \ref{s:TPCA} we study the spiked tensor model, as an application of our results on low rank perturbations of long random matrices. The proofs of our main results are given in Section \ref{s:mainlong}.

\noindent\textbf{Notations }
For two numbers $X,Y$, we write that $X = \OO(Y )$ if there exists a universal constant $C>0$ such
that $|X| \leq C Y$. We write $X = \oo(Y )$, or $X \ll Y$ if the ratio $|X|/Y\rightarrow \infty$ as $n$ goes to infinity. We write
$X\asymp Y$ if there exists a universal constant $C>0$ such that $Y/C \leq |X| \leq  C Y$.
We denote the index set $\qq{k}=\{1,2,\cdots,k\}$.
 We say an event $\Omega$ holds with high probability if for any large $C>0$, $\bP(\Omega)\geq 1-n^{-C}$ for $n$ large enough. 
We write $X\prec Y$ or $X=\OO_\prec(Y)$, if $X$ is stochastically dominated by $Y$ in the sense that for all large $C>0$, we have
\begin{align*}
\bP(X\geq n^{1/C}Y)\leq n^{-C},
\end{align*}
for large enough $n$.

\noindent\textbf{Acknowledgements }
The research of J.H. is supported by the Simons Foundation as a Junior Fellow at
the Simons Society of Fellows, and NSF grant DMS-2054835.

\section{Main Results}\label{s:main}

In this section, we state our main results on the singular values and singular vectors of low rank perturbations of long random matrices. Let $Z\in \bR^{n\times m}$ be an $n\times m$ random matrices, with i.i.d. random entries $Z_{ij}$ satisfying
\begin{assumption}\label{a:Zasup}
The entries of $Z$ are i.i.d., and have mean zero and variance $1/n$, and for any integer $p\geq 1$,
\begin{align*}
\bE[Z_{ij}]=0, \quad \bE[Z_{ij}^2]=\frac{1}{n}, \quad \bE[Z_{ij}^{2p}]\leq \frac{C_p}{n^p},\quad 1\leq i\leq n,\quad 1\leq j\leq m.
\end{align*}
\end{assumption}
We introduce the parameter $\phi:=\sqrt{m/n}$.
It is well-known that for $m,n\rightarrow\infty$, if their ratio converges to $m/n=\phi^2\rightarrow \phi_*^2\in[1,\infty)$ which is independent of $n$, the empirical eigenvalue distribution of $Z Z^\top$ converges to the Marchenko-Pastur law
\begin{align}\label{e:MPa}
\frac{\sqrt{\left(x-\left(\phi_*-1\right)^2\right)\left(\left(\phi_*+1\right)^2-x\right)}}{2\pi  x}\rd x,
\end{align}
supported on $[(\phi_*-1)^2, (\phi_*+1)^2]$.

In this paper, we consider the case where the ratio $m/n=\phi^2$ depends on $n$, satisfying that $n^{1/C}\leq \phi\leq n^C$ for some constant $C>0$. In Section \ref{s:localZ}, we collect some results on singular values of the long random matrix $Z$ from \cite{alex2014isotropic}.
Using the estimates of the singular values and singular vectors of long random matrices as input, we can study low rank perturbations of long random matrices. We consider  rank one perturbation of long random matrices in the form,
\begin{align}\label{e:bvuc}
\beta \bmv \bmu^\top+ Z\in \bR^{n\times m},
\end{align}
where $\bmv\in \bR^n$, $\bmu\in \bR^m$ are unit vectors, and $Z$ is an $n\times m$ random matrix satisfying Assumption \ref{a:Zasup}.
We state our main results on the singular values and singular vectors of low rank perturbations of long random matrices \eqref{e:bvuc} in Section \ref{s:lowrank}.

\subsection{Singular values of Long Random Matrices}\label{s:localZ}
Let $Z$ be an $n\times m$ random matrix, with entries satisfying Assumption \ref{a:Zasup}.  Let $\phi=\sqrt{m/n}$ satisfy that $n^{1/C}\leq \phi\leq n^C$ for some constant $C>0$. The eigenvalues of sample covariance matrices $ZZ^\top$ in this setting have been well studied in \cite{alex2014isotropic}. We denote the singular values of $Z$ as
\begin{align*}
s_1\geq s_2\geq \cdots \geq s_n\geq 0.
\end{align*}
They are square roots of eigenvalues of the sample covariance matrices $ZZ^\top$. As an easy consequence of \cite[Theorem 2.10]{alex2014isotropic}, we have the following theorem on the estimates of the largest singular value of $Z$.

\begin{theorem}\label{t:extreme}
Under Assumption \ref{a:Zasup}, let $\phi=\sqrt{m/n}$ with  $n^{1/C}\leq \phi\leq n^C$. Fix an arbitrarily small $\delta>0$, with high probability the largest singular value $s_1$ of $Z$ satisfies
\begin{align}\label{e:s1est}
|s_1-(\phi+1)|\leq \frac{n^\delta}{n^{2/3}},
\end{align}
provided $n$ is large enough.
\end{theorem}

\begin{figure}[H]
\center
\includegraphics[scale=0.6]{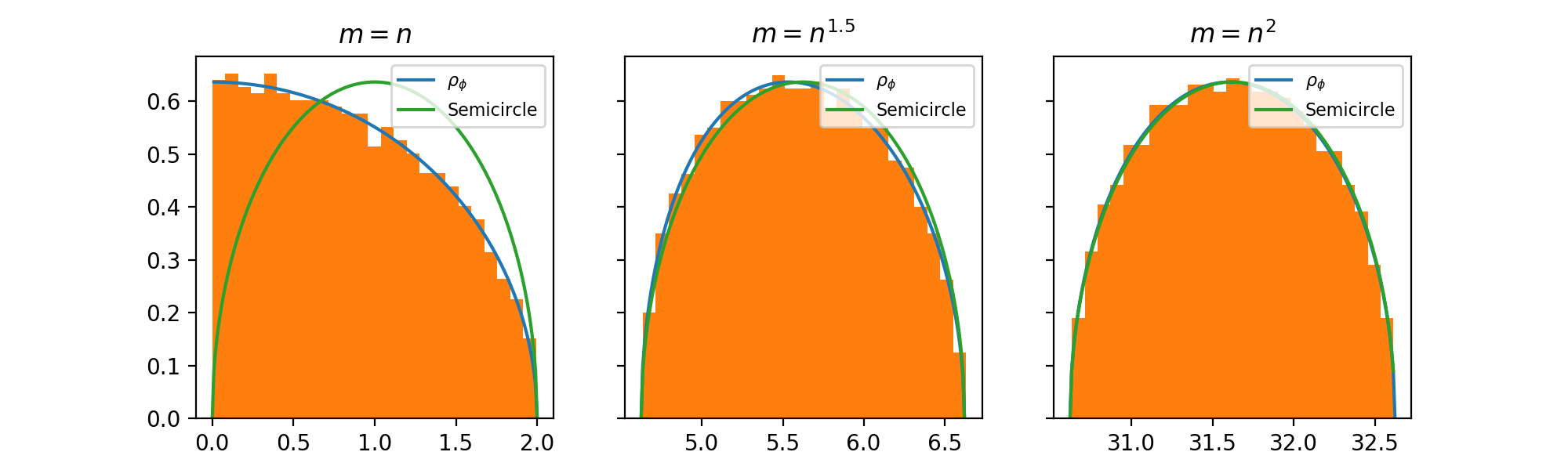}
\caption{The empirical singular value distribution of $Z$, with $n=1000$ and $m=n, n^{1.5}, n^2$. When $\phi=\sqrt{m/n}$ is large, the empirical singular value converges to the Semicircle distribution supported on $[\phi-1, \phi+1]$.} \label{f:sqMP}
\end{figure}

The results in \cite{alex2014isotropic} give estimates of each eigenvalues of the sample covariance matrix $ZZ^\top$ away from $0$, see Theorem \ref{t:rigidity}. It also gives estimates of locations of each singular value of $Z$. In particular, it implies that the empirical singular value distribution $\sum_{i=1}^n \delta_{s_i}/n$ of $Z$ is close to the pushforward of the Marchenko-Pastur law (after proper normalization)  by the map $x\mapsto \sqrt{x}$,  
\begin{align*}
&\phantom{{}={}}\rho_{\phi}(x)\rd x
=\frac{\sqrt{\left(x^2-\left(\phi-1\right)^2\right)\left(\left( \phi+1\right)^2-x^2\right)}}{\pi x}\rd x.
\end{align*}
We remark that $\rho_\phi$ depends on $m,n$ through $\phi$ and is supported on $[\phi-1, \phi+1]$. As $\phi\rightarrow \infty$ with $n$, $\rho_\phi(x)$ after shifting by $\phi$, i.e. $\rho_\phi(x+\phi)$, converges to the semi-circle distribution on $[-1, +1]$:
\begin{align*}
\rho_\phi(x+\phi)\rd x\rightarrow \frac{2\sqrt{(x+1)(1-x)}}{\pi}.
\end{align*}
See Figure \ref{f:sqMP} for some plots of $\rho_\phi$. One can see that the extreme singular values stick to the boundary of the support of the limiting empirical measure.

\subsection{Low Rank Perturbations of Long Random Matrices}\label{s:lowrank}
Let $Z$ be an $n\times m$ random matrix, with entries satisfying Assumption \ref{a:Zasup}.  Let $\phi=\sqrt{m/n}$. Without loss of generality we can assume that $\phi\geq 1$, otherwise, we can simply study the transpose of $Z$. We allow $\phi$ to grow with $n$ at any polynomial rate, $1\leq \phi\leq n^C$ for any large number $C>0$. In this regime, $Z$ is a long random matrix. 

In this section, we state our main results on the rank one perturbation of long random matrices from \eqref{e:bvuc}:
\begin{align}\label{e:bvuccopy}
\beta \bmv \bmu^\top+ Z\in \bR^{n\times m},
\end{align}
where $\bmv\in \bR^n$, $\bmu\in \bR^m$ are unit vectors.
As in Theorem \ref{t:extreme}, in this setting, the singular values of $Z$ are roughly supported on 
the interval $[\phi-1, \phi+1]$.
The following theorem states that there is an exact $n$-dependent threshold $\sqrt{\phi}$, if $\beta$ is above the threshold, $\beta\bmv\bmu^\top+Z$ has an outlier singular value; if $\beta$ is below this threshold, there are no outlier singular values, and all the singular values are stick to the bulk.
\begin{theorem}\label{t:eig}
We assume Assumption \ref{a:Zasup} and $1\leq \phi=\sqrt{m/n}\leq n^C$. Let $\beta=\la\sqrt\phi$ with $\la\geq 0$, fix arbitrarily small $\fc>0$, and denote $\hat s_1$ the largest singular value of $\beta \bmv \bmu^\top + Z$. For any small $\delta>0$, if $\la\geq 1+n^{-1/3+\fc}$, with high probability, the largest singular value $\hat s_1$ of $\beta \bmv \bmu^\top + Z$ is an outlier, and explicitly given by
\begin{align}\label{e:strong}
\hat s_1=\sqrt{\phi^2+(\la^2+1/\la^2)\phi+1}+\OO \left(\frac{n^\delta (\la-1)^{1/2}}{n^{1/2} \phi}\right).
\end{align}
If $\la\leq 1+n^{-1/3+\fc}$, with high probability, $\beta \bmv \bmu^\top + Z$ does not have outlier singular values, and the largest singular value $\hat s_1$ satisfies
\begin{align}\label{e:weak}
\hat s_1\leq \phi+1+n^{-2/3+3\fc},
\end{align}
provided $n$ is large enough.
\end{theorem}

\begin{figure}[H]
\center
\includegraphics[scale=0.6]{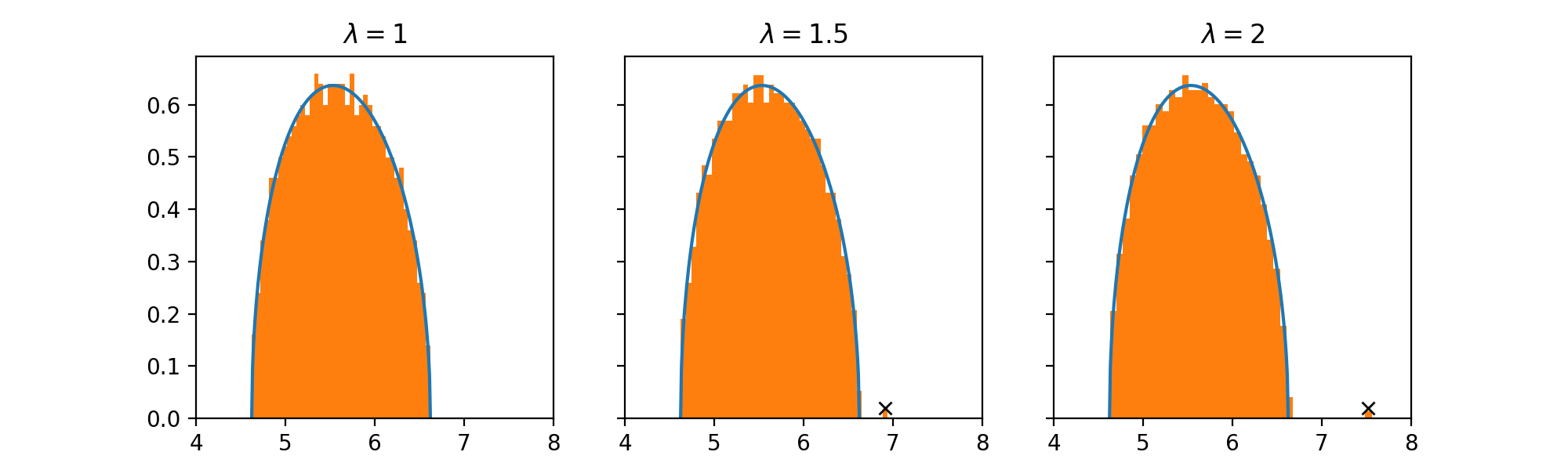}
\caption{The empirical singular value distribution of $\beta \bmv\bmu^\top+Z$, with $n=1000$,$m=n^{1.5}$, and $\beta=\lambda (m/n)^{1/4}$ for $\lambda=1, 1.5, 2$. We marked the predicted outlier by $\times$ as given by formula  \eqref{e:strong}.} \label{f:Outlier}
\end{figure}

We refer to Figure \ref{f:Outlier} for an illustration of Theorem \ref{t:eig}. Theorem \ref{t:eig} also characterizes the behavior of the outlier in the critical case, when $\lambda$ is close to $1$. 

We have similar transition for the singular vectors. If $\beta$ is above the threshold $\sqrt{\phi}$, the left singular vector $\hat\bmv_1$ associated with the largest singular value of $\beta \bmv \bmu^\top + Z$ has a large component in the signal $\bmv$ direction; If $\beta$ is below the threshold $\sqrt{\phi}$, the projection of  $\hat\bmv_1$ on the signal $\bmv$ direction vanishes. 
\begin{theorem}\label{t:ev}
We assume Assumption \ref{a:Zasup} and $1\leq \phi=\sqrt{m/n}\leq n^C$. 
Let $\beta=\la\sqrt\phi$ with $\la\geq 0$, fix arbitrarily small $\fc>0$, and denote $\hat s_1$ the largest singular value of $\beta \bmv \bmu^\top + Z$. For any small $\delta>0$, if $\la\geq 1+n^{-1/3+\fc}$, with high probability, the left singular vector $\hat\bmv_1$ associated with the largest singular value of $\beta \bmv \bmu^\top + Z$ has a large component in  the signal $\bmv$ direction:
\begin{align}\label{e:strong1}
|\langle \hat\bmv_1, \bmv\rangle |=\left(1+\OO\left(\frac{n^\delta}{n^{1/2}(\la-1)^{3/2}}\right)\right)\sqrt{\frac{(\la^4-1)}{\la^4+\la^2/\phi}}.
\end{align}
And similar estimates for the right singular vector $\hat\bmu_1$ associated with the largest singular value of $\beta \bmv \bmu^\top + Z$
\begin{align}\label{e:strong22}
&|\langle \hat{\bmu}_1, \bmu\rangle| 
=\left(1+\OO_\prec\left(\frac{1}{n^{1/2}(\la-1)^{3/2}}\right)\right)\sqrt{\frac{\la^4-1}{\la^2(\la^2+\phi)}}.
\end{align}
If $\la\leq 1+n^{-1/3+\fc}$, with high probability, the projection of $\hat \bmv_1$ on $\bmv$, and the projection of $\hat\bmu_1$ on $\bmu$ are upper bounded by
\begin{align}\label{e:weak1}
\max\{| \langle \hat \bmv_1, \bmv\rangle|,|\langle \hat \bmu_1, \bmu\rangle|\}\leq  n^{4\fc}\min\left\{n^{-1/6}, \frac{1}{\sqrt n(\la-1)}\right\},
\end{align}
provided $n$ is large enough. 
\end{theorem}

\begin{remark}
We remark that when the ratio $\phi=\sqrt{m/n}$ goes to infinity with $n$, with $\lambda$ fixed, then from \eqref{e:strong1}, $|\langle \hat \bmv_1, \bmv\rangle|$ converges to $\sqrt{1-\la^{-4}}$, and from \eqref{e:strong22}, $|\langle \hat \bmu_1, \bmu\rangle|$ converges to $0$.
\end{remark}

\begin{remark}
In this paper, for simplicity of notations, we only consider rank one perturbations of long random matrices. Our method can as well be used to study any finite rank perturbations of long random matrices.
\end{remark}

The singular values and vectors of low rank perturbations of large rectangular random matrices has previously been studied in \cite{benaych2012singular} by Benaych-Georges and Nadakuditi. Our main results Theorems \ref{t:eig} and \ref{t:ev} are generalization of their results in two directions. Firstly, in \cite{benaych2012singular},  the ratio $\phi=\sqrt{m/n}$ is assumed to converge to a constant independent of $n$. In our setting, we allow $\phi$ to have polynomial growth in $n$. As we will see in Section \ref{s:TPCA}, this will be crucial for us to study the tensor principle component analysis. Secondly, we allow the signal-to-noise ratio to be close to the threshold (depending on $n$), and our main results characterize the behaviors of singular values and singular vectors in this regime.  In an independent work \cite{feldman2021spiked} by Feldman, the singular values and vectors of multi-rank perturbations of long random matrices have been studied under the assumption that either the signal vectors contains i.i.d. entries with mean zero and variance one, or the noise matrix has Gaussian entries. Both proofs use the master equations which characterize the singular values and singular vectors of low rank perturbations of rectangular random matrices developed in \cite{benaych2012singular}, see Section \ref{s:master}. To analyze the master equation, in \cite{feldman2021spiked}, Feldman needs that the signal vectors have i.i.d. entries with mean zero and variance one. Our argument uses results  from \cite{alex2014isotropic}, which gives Green's function estimates of long random matrices, and works for any (deterministic) signal vectors.

\section{Tensor PCA}\label{s:TPCA}
As an application of our main results Theorems \ref{t:eig} and \ref{t:ev}, in this section, we use them to study the asymmetric rank-one spiked tensor model as introduced in \cite{NIPS2014_5616}:
\begin{align}\label{e:rank1}
\bmX=\beta \bmv_1\otimes \bmv_2\otimes \cdots \otimes \bmv_k+\bmW,
\end{align}
where
\begin{itemize}
\item $\bmX\in \otimes^k \bR^n$ is the $k$-th order tensor observation.
\item $\bmW\in \otimes^k \bR^n$ is a noise tensor. The entries of $\bmW$ are independent random variables with mean zero and variance $1/n$.
\item $\beta\in \bR$ is the signal-to-noise ratio.
\item $\bmv_i\in \bR^n$ are unknown unit vectors to be recovered.
\end{itemize}
The goal is 
to perform reliable estimation and inference on the unseen signal tensor $\bmv_1\otimes \bmv_2\otimes \cdots \otimes \bmv_k$. We remark that for any rank-one tensor $T\in \otimes^k\bR^n$, it can be uniquely written as $T=\beta \bmv_1\otimes \bmv_2\otimes \cdots \otimes \bmv_k$, where $\bmv_i$ are unit vectors. The model \eqref{e:rank1} is slightly more general than the asymmetric spiked tensor model in \cite{NIPS2014_5616}, which assumes that $\bmv_1=\bmv_2=\cdots=\bmv_k$. In this section, we make the following assumption on the noise tensor:
\begin{assumption}\label{a:Wasup}
The entries of $\bmW$ are i.i.d., and they have mean zero and variance $1/n$: for any indices $1\leq i_1,i_2,\cdots, i_k\leq n$, and any integer $p\geq 1$,
\begin{align*}
\bE[\bmW_{i_1i_2\cdots i_k}]=0, \quad \bE[\bmW_{i_1i_2\cdots i_k}^2]=\frac{1}{n}, \quad \bE[\bmW_{i_1i_2\cdots i_k}^{2p}]\leq \frac{C_p}{n^p}.
\end{align*}
\end{assumption}

The tensor unfolding algorithm in  \cite{NIPS2014_5616} unfolds the tensor $\bmX$ to an $n^{\floor{k/2}}\times n^{\ceil{k/2}} $ matrix, and they proved that it detects the signal when the signal-to-noise ratio satisfies $\beta\gg n^{(\ceil{k/2}-1)/2}$. The conjectured algorithmic threshold is $\beta\gg n^{(k-2)/4}$. We recall the tensor unfolding algorithm. Take any index set $\bI\subseteq \qq{1,k}$, with $1\leq |\bI|=q\leq k/2$. Given any $\bI=\{k_1,k_2,\cdots, k_q\}$, let $\qq{1,k}\setminus \bI=\{\ell_1, \ell_2,\cdots, \ell_{k-q}\}$, we denote the matrix ${\Mat_\bI(\bmX)}$, which is the $n^q\times n^{k-q}$ matrix obtained from $\bmX$ by unfolding along the axes indexed by $\bI$. More precisely, for any indices $i_1,i_2,\cdots, i_k\in \qq{1,n}$, let
$a=1+\sum_{j=1}^q(i_{k_j}-1)n^{j-1}$ and $b=1+\sum_{j=1}^{k-q}(i_{\ell_j}-1)n^{j-1}$, then
\begin{align*}
\Mat_\bI(\bmX)_{a,b}=\bmX_{i_1i_2\cdots i_k}.
\end{align*}
If $\bI$ is a singleton, i.e. $\bI=\{i\}$, we will simply write $\Mat_\bI$ as $\Mat_i$.

We can view $\Mat_\bI(\bmX)$ as the sum of the unfolding of the signal tensor $ \bmv_1\otimes \bmv_2\otimes \cdots \otimes \bmv_k$ and the noise tensor $\bmW$. Let
\begin{align*}
\bmv_\bI\deq \text{vec}[\otimes_{j=1}^q \bmv_{k_j}]\in \bR^{n^q},\quad \bmu_{\bI}\deq \text{vec}[\otimes_{j=1}^{k-q} \bmv_{\ell_j}]\in  \bR^{n^{k-q}},\quad
Z_\bI\deq \Mat_\bI(\bmW)\in  \bR^{n^q}\times \bR^{n^{k-q}}.
\end{align*}
Then we can rewrite $\Mat_\bI(\bmX)$ as 
\begin{align}\label{e:Mat}
\Mat_\bI(\bmX)=\beta \bmv_\bI \bmu^\top_{\bI}+Z_\bI(\bmW)\in  \bR^{n^q}\times \bR^{n^{k-q}}.
\end{align}
To make \eqref{e:Mat} in the form of \eqref{e:bvuc}, we need to further normalize $\Mat_\bI(\bmX)$ as
\begin{align}\label{e:newMat}
\widetilde{\Mat_\bI}(\bmX)=\frac{\beta}{n^{(q-1)/2}} \bmv_\bI \bmu^\top_{\bI}+\widetilde{Z_\bI},\quad \widetilde{Z_\bI}=\frac{{Z_\bI}}{n^{(q-1)/2}}.
\end{align}
In this way, each entry of $\widetilde{Z_\bI}$ has variance $1/n^{q}$. And \eqref{e:newMat} is a rank one perturbation of a $n^{q}\times n^{k-q}$ random matrix $\widetilde{Z_\bI}$ in the form of \eqref{e:bvuc} (by taking $(n,m)$ as $(n^q, n^{k-q})$), and the ratio $\phi=\sqrt{n^{k-q}/n^q}=n^{(k-2q)/2}\geq 1$ grows at most polynomially in $n$.  We take 
$\hat s_1$ the largest singular value of the normalized unfolded matrix $\widetilde{\Mat_\bI}(\bmX)$, our main results Theorems \ref{t:eig} and \ref{t:ev} indicate that there is a phase transition at $\beta=n^{(k-2)/4}$ for the tensor unfolding \eqref{e:newMat}.

\begin{theorem}\label{t:tensor}We assume Assumption \ref{a:Wasup}, and fix any index set $\bI\subseteq \qq{1,k}$ with $|\bI|=q\leq k/2$. Let 
$\widetilde{\Mat_\bI}(\bmX)$ be the normalized matrix obtained from $\bmX$ by unfolding along the axes indexed by $\bI$, as in \eqref{e:newMat}, and denote the ratio $\phi=\sqrt{n^{k-q}/n^q}=n^{(k-2q)/2}\geq 1$.

Let $\beta=\la n^{(k-2)/4}$ with $\la\geq 0$, fix arbitrarily small $\fc>0$, and denote $\hat s_1$ the largest singular value of $\widetilde{\Mat_\bI}(\bmX)$, and $\hat\bmv_1$ the corresponding left singular vector. For arbitrarily small $\delta>0$, if $\la\geq 1+n^{-(1/3-\fc)q}$, with high probability, the largest singular value $\hat s_1$ is an outlier, is explicitly given by
\begin{align*}
\hat s_1=\sqrt{\phi^2+(\la^2+1/\la^2)\phi+1}+\OO \left(\frac{n^\delta (\la-1)^{1/2}}{n^{q/2} \phi}\right),
\end{align*}
and the left singular vector $\hat\bmv_1$ has a large component in the $\bmv_\bI$ direction:
\begin{align}\label{e:strong13}
|\langle \hat\bmv_1, \bmv_{\bI}\rangle |=\left(1+\OO\left(\frac{n^\delta}{n^{q/2}(\la-1)^{3/2}}\right)\right)\sqrt{\frac{(\la^4-1)}{\la^4+\la^2/\phi}}.
\end{align}

If $\la\leq 1+n^{-(1/3-\fc)q}$, with high probability, $\widetilde{\Mat_\bI}(\bmX)$ does not have outlier singular values, the largest singular value $\hat s_1$ satisfies
\begin{align*}
\hat s_1\leq \phi+1+n^{-(2/3-3\fc)q},
\end{align*}
and the projection of $\hat \bmv_1$ on $\bmv_\bI$ is upper bounded by
\begin{align}\label{e:strong23}
|\langle \hat\bmv_1, \bmv_\bI\rangle|\leq  n^{4q\fc}\min\left\{n^{-q/6}, \frac{1}{n^{q/2}(\la-1)}\right\},
\end{align}
provided $n$ is large enough.
\end{theorem}
\begin{proof}[Proof of Theorem \ref{t:tensor}]
Theorem \ref{t:tensor} follows from Theorems \ref{t:eig} and \ref{t:eig} by taking $(n,m)$ as $(n^q, n^{k-q})$, and the signal-to-noise ratio $\beta$ as $\beta/n^{(q-1)/2}$. In this way the criteria in Theorems \ref{t:eig} and \ref{t:eig} become that if 
\begin{align*}
\frac{\beta/n^{(q-1)/2}}{\sqrt\phi}=\frac{\beta}{n^{(q-1)/2}n^{(k-2q)/4}}=\frac{\beta}{n^{(k-2)/4}}=\la\geq 1+n^{-(1/3-\fc)q}
\end{align*} 
then $\hat s_1$ is an outlier with high probability and the left singular vector $\hat\bmv_1$ has a large component in the $\bmv_\bI$ direction. Otherwise if 
\begin{align*}
\la\leq 1+n^{-(1/3-\fc)q},
\end{align*}
$\hat s_1$ sticks to the bulk and the projection of $\hat \bmv_1$ on $\bmv_\bI$ is small.
\end{proof}

We remark that as indicated by Theorem \ref{t:tensor}, the tensor unfolding algorithms which unfold $\bmX$ to an $n^q\times n^{k-q}$ matrix for any choice of $1\leq q\leq k/2$ and index set $\bI$,  essentially share the same threshold, i.e. $\beta=n^{(n-2)/4}$, which matches the conjectured threshold. 
As in \eqref{e:strong13}, above the signal-to-noise ratio threshold, the left singular vector $\hat\bmv_1$  corresponding to the largest singular value of $\widetilde{\Mat_\bI}(\bmX)$ is aligned with $\bmv_\bI$. The leading order accuracy for the estimator $\hat \bmv_1$ as in \eqref{e:strong13} is the same and is independent of $q$.
However, if $q>1$, this does not give us information of individual signal $\bmv_i$. A further recursive unfolding method is proposed in \cite{NIPS2014_5616} to recover each individual signal $\bmv_i$.

Since by taking $q>1$ does change the algorithmic threshold, but increasing the computation cost. 
We propose the following simple algorithm to recover each signal $\bmv_i$, by performing the tensor unfolding algorithm for each $1\leq i\leq k$ with $\bI=\{i\}$, namely $q=1$.

For each $1\leq i\leq k$, we unfold \eqref{e:rank1} to an $n\times n^{k-1}$ matrix
\begin{align}\label{e:unfold1}\begin{split}
&\widetilde{\Mat_{i}}(\bmX)={\Mat_{i}}(\bmX)=\beta \bmv_i \bmu^\top_i + Z_i,\\
& \bmu_i=\text{vec}[\bmv_1\otimes \bmv_2\otimes \cdots \otimes \bmv_{i-1}\otimes \bmv_{i+1}\otimes\cdots \otimes \bmv_k]\in \bR^{n^{k-1}},
\end{split}\end{align}
which is \eqref{e:newMat} by taking $q=1$ and $\bI=\{i\}$.
In this way \eqref{e:unfold1} is a rank one perturbation of a long $n\times n^{k-1}$  random matrix in the form of \eqref{e:bvuc}, and the ratio $\phi=\sqrt{n^{k-1}/n}=n^{(k-1)/2}$ grows with $n$. We take 
$\hat s^{(i)}_1$ the largest singular value of ${\Mat_{i}}(\bmX)=\beta \bmv_i\otimes \bmu_i^\top + Z_i$, and denote
 \begin{align}\label{e:defhbb}
 \hat\beta^{(i)}=\sqrt{\frac{((\hat s^{(i)}_1)^2-(\phi^2+1))+\sqrt{((\hat s^{(i)}_1)^2-(\phi+1)^2)((\hat s^{(i)}_1)^2-(\phi-1)^2)}}{2}},
 \end{align}
 as the estimator for $\beta$; and 
$\hat\bmv^{(i)}_1$ the left singular vector corresponding to the largest singular value of ${\Mat_{i}}(\bmX)=\beta \bmv_i \bmu_i^\top + Z_i$, as the estimator for $\bmv_i$. This gives the following simple algorithm to recover $\beta$ and $\bmv_i$.

\begin{algorithm}[H]
\SetAlgoLined
\textbf{Input:} $\bmX$\;
 \For{$i$ from $1$ to $k$}{
 $\hat s_1^{(i)}=$ largest singular value of $\Mat_i(\bmX)$\;
 $\bmv_1^{(i)}=$ left singular vector corresponding to the largest singular value of $\Mat_i(\bmX)$\;
 $\hat\beta^{(i)}=\sqrt{\left(((\hat s^{(i)}_1)^2-(\phi^2+1))+\sqrt{((\hat s^{(i)}_1)^2-(\phi+1)^2)((\hat s^{(i)}_1)^2-(\phi-1)^2)}\right)/2}$\;
 }
   \KwResult{$\{\hat\beta^{(i)}, \bmv_1^{(i)}\}_{1\leq i\leq k}$.}

 \caption{Tensor Unfolding}
\end{algorithm}
In tensor literature, the above algorithm is exactly the truncated higher order singular value decomposition (HOSVD) introduced in \cite{de2000multilinear}. The higher order orthogonal iteration (HOOI), which uses the truncated HOSVD as initialization combining with a power iteration, was developed in \cite{de2000best} to find a best low-multilinear-rank approximation of a tensor. The performance of HOOI was analyzed in \cite{zhang2018tensor} for the spiked tensor model. It was proven that for the signal-to-noise ratio $\beta=\la n^{(k-2)/4}$ with $\la\geq C_{gap}$ for some large constant $C_{gap}>0$, HOOI converges within a logarithm factor of iterations. As an easy consequence of Theorem \ref{t:tensor}, we have the following theorem which gives the exact threshold of the signal-to-noise ratio, i.e. $\beta=n^{(n-2)/2}$. Above the threshold, our estimators $\hat\beta^{(i)}$ and $\hat\bmv_1^{(i)}$ approximate the signal-to-noise ratio $\beta$ and the signal vector $\bmv_i$. 
\begin{theorem}\label{t:tensor2}
We assume Assumption \ref{a:Wasup}. For any $1\leq i\leq k$, we unfold $\bmX$ to $\beta \bmv_i \bmu^\top_i + Z_i$ as in \eqref{e:unfold1}.
Let the estimator $\hat\beta^{(i)}$ be as defined in \eqref{e:defhbb}, and $\hat\bmv^{(i)}_1$ the left singular vector corresponding to the largest singular value of ${\Mat_{i}}(\bmX)=\beta \bmv_i \bmu_i^\top + Z_i$.

Let $\beta=\la n^{(k-2)/4}$ with $\la\geq 0$, and fix arbitrarily small $\fc>0$. For arbitrarily small $\delta>0$, if $\la\geq 1+n^{-1/3+\fc}$, with high probability, $\hat\beta^{(i)}$ and $\hat\bmv_1^{(i)}$ approximates $\beta$ and $\bmv_i$
\begin{align}\label{e:strongc}
|\hat\beta^{(i)}-\beta|\leq  \frac{n^\delta}{n^{(k+1)/4}\sqrt{\la-1}},
\end{align}
and 
\begin{align}\label{e:strong1c}
|\langle \hat\bmv^{(i)}_1, \bmv_i\rangle| =\left(1+\OO\left(\frac{n^\delta}{n^{1/2}(\la-1)^{3/2}}\right)\right)\sqrt{\frac{(\la^4-1)}{\la^4+\la^2 n^{-(k-1)/2}}}
\sim \sqrt{1-\frac{1}{\la^4}}.
\end{align}
If $\la\leq1+n^{-1/3+\fc}$, with high probability, the projection of $\hat \bmv_1^{(i)}$ on $\bmv_i$ vanishes as $\la$ decreases
\begin{align}\label{e:weak1c}
|\langle \hat\bmv_1^{(i)}, \bmv_i\rangle|\leq n^{4\fc} \min\left\{n^{-1/6}, \frac{1}{\sqrt n|\la-1|}\right\},
\end{align}
provided $n$ is large enough.
\end{theorem}

\begin{remark}
Given the unfolded $n\times n^{k-1}$ matrix ${\Mat_{i}}(\bmX)=\beta \bmv_i\otimes \bmu_i^\top + Z_i$, the largest singular value and its left eigenvector can be computed by first computing ${\Mat_{i}}(\bmX){\Mat_{i}}(\bmX)^\top$, then computing the the largest eigenvalue and corresponding eigenvector by power iteration. The total time complexity is $\OO(Tn^k)$, where $T$ is the number of iterations, which can be taken as $T=\ln n$. Therefore, the estimators $\hat \beta^{(i)}$ and $\hat \bmv_1^{(i)}$ can be computed  with time complexity $\OO((\ln n) n^k)$.
To recover the signals $\bmv_i$ for each $1\leq i\leq k$, we need to repeat the above tensor unfolding algorithm $k$ times, and obtain $\bmv_1^{(i)}$ for each $1\leq i\leq k$. The total time complexity is $\OO((\ln n) kn^k)$.
\end{remark}

\begin{proof}
The claims \eqref{e:strong1c} and \eqref{e:weak1c} follow directly from \eqref{e:strong13} and \eqref{e:strong23} by taking $q=1$ and $\phi=n^{(k-1)/2}$. In the following we prove \eqref{e:strongc}. Fix arbitrarily small $\delta>0$. For $\la\geq 1+n^{-1/3+\fc}$, with high probability, the largest singular value $\hat s_1^{(i)}$ of ${\Mat_{i}}(\bmX)=\beta \bmv_i \bmu_i^\top + Z_i$ is given by \eqref{e:strong}
\begin{align}\begin{split}\label{e:diffa}
\hat s_1^{(i)}&=\sqrt{\phi^2+(\la^2+1/\la^2)\phi+1}+\OO \left(\frac{n^\delta(\la-1)^{1/2}}{n^{1/2} \phi}\right)\\
&=\sqrt{\phi^2+\beta^2+\phi^2/\beta^2+1}+\OO\left(\frac{n^\delta(\la-1)^{1/2}}{n^{1/2} \phi}\right),
\end{split}\end{align}
where $\phi=n^{(k-2)/4}$. Our $\hat \beta^{(i)}$ as defined in \eqref{e:defhbb} is chosen as the solution of 
\begin{align}\label{e:diffhi}
\hat s^{(i)}_1&=\sqrt{\phi^2+ (\hat\beta^{(i)})^2+\phi^2/(\hat\beta^{(i)})^2+1}.
\end{align}
By taking difference of \eqref{e:diffa} and \eqref{e:diffhi}, and rearranging, we get
\begin{align}\label{e:difhihi}
\beta^2+\frac{\phi^2}{\beta^2}-(\hat\beta^{(i)})^2-\frac{\phi^2}{(\hat\beta^{(i)})^2}
=(\beta^2-(\hat\beta^{(i)})^2)\left(1-\frac{\phi^2}{\beta^2(\hat\beta^{(i)})^2}\right)
\leq \frac{n^\delta (\la-1)^{1/2}}{n^{1/2}}\frac{\phi+\beta}{\phi}.
\end{align}
We notice that $\beta^2(\hat\beta^{(i)})^2/\phi^2-1\gtrsim (\la^4-1)\gtrsim (\la-1)$. Thus \eqref{e:difhihi} implies
\begin{align*}
|\beta-(\hat\beta^{(i)})|\leq \frac{n^\delta}{n^{1/2}\sqrt \phi \sqrt{\la-1}},
\end{align*}
with high probability, provided $n$ is large enough.
This finishes the proof \eqref{e:strongc}.
\end{proof}

\begin{figure}
\includegraphics[scale=0.6]{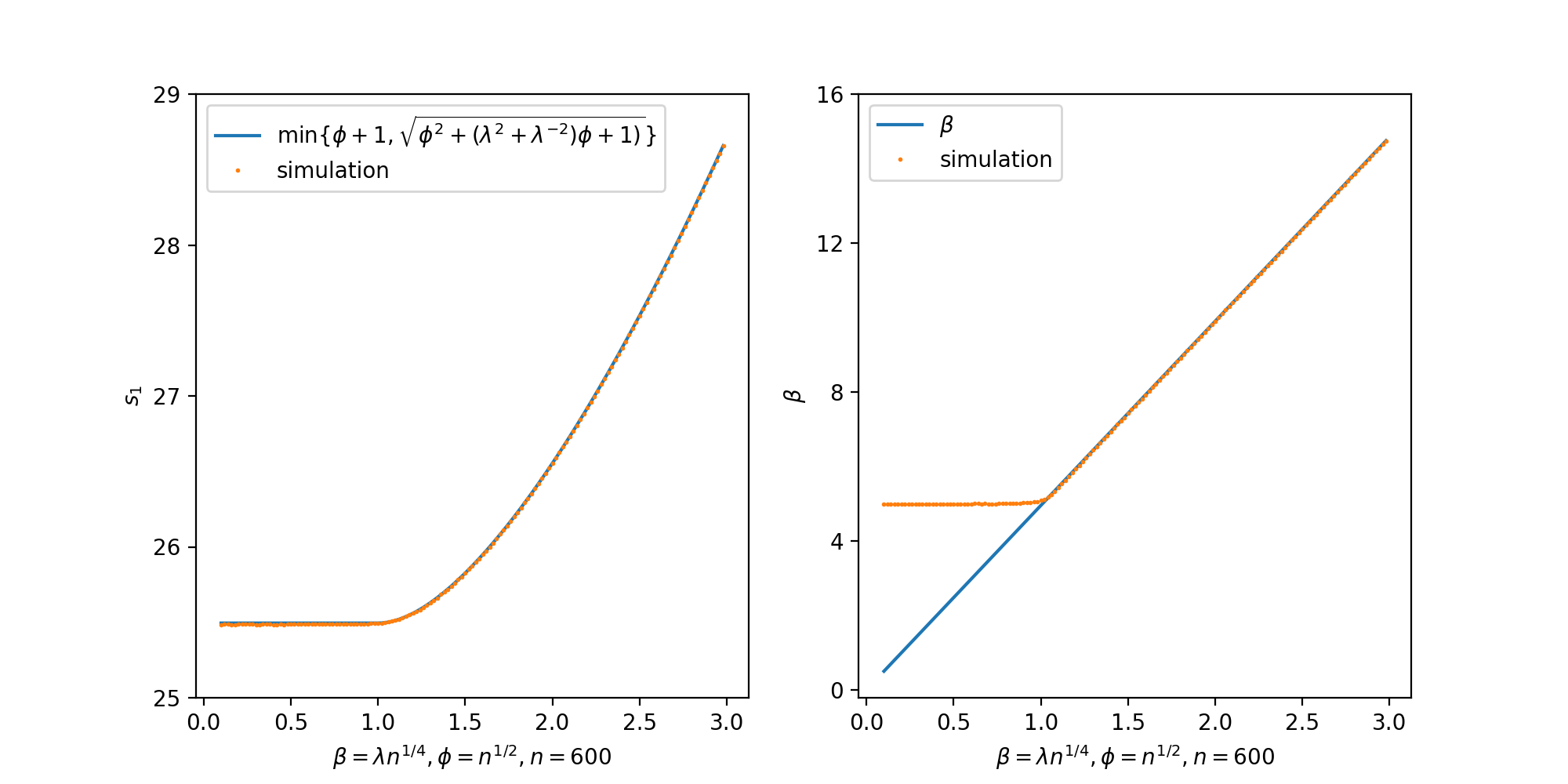}
\caption{For $n=600$, $k=3$, $\beta=\la*n^{(k-2)/4}$ with $\la\in[0,3]$, averaging over $500$ trials, the left panel plots the largest singular value of the unfolded matrix and our theoretical prediction; the left panel plots the $\beta$ and the estimator $\hat\beta$.} \label{f:eig_plot}
\end{figure}

\begin{figure}
\includegraphics[scale=0.6]{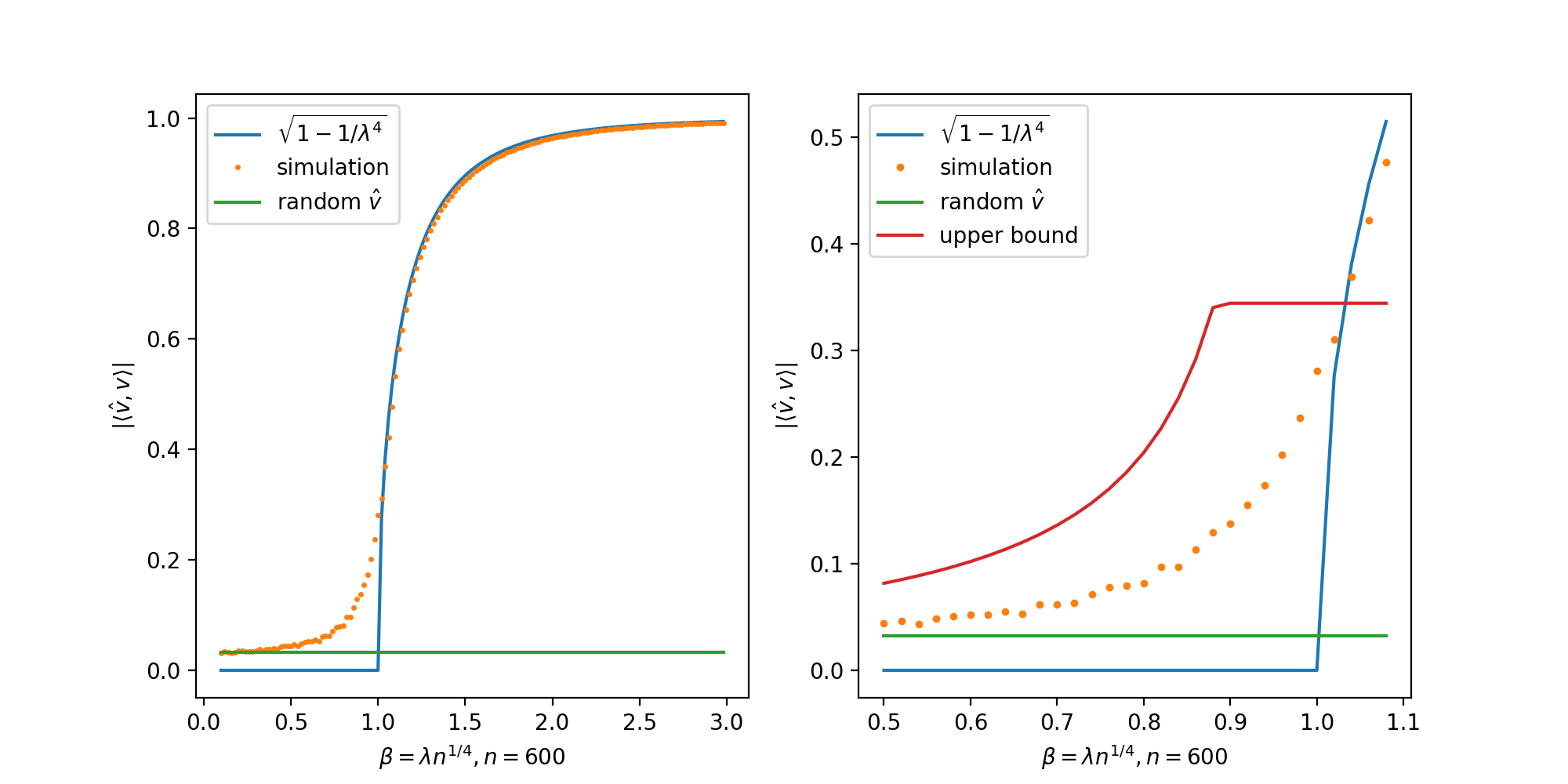}
\caption{This plot of $|\langle \hat{\bmv}, \bmv\rangle|$ is for $n=600$, $\beta=\la*n^{(k-2)/4}$ with $\la\in[0,3]$, averaging over $500$ trials.} \label{f:ev_plot}
\end{figure}

We numerically verify Theorem \ref{t:tensor}. We take $n=600$, $k=3$ and $\beta=\la n^{(k-2)/4}=\la n^{1/4}$
for $\la\in[0,3]$. We sample the signals $\bmv_1=\bmv_2=\bmv_3=\bmv$ as unit Gaussian vectors, and the noise tensor $\bmW$ with independent Gaussian entries. In the left panel of Figure \ref{f:eig_plot}, we plot the largest singular value of the unfolded matrix ${\Mat_{i}}(\bmX)$ and our theoretical prediction \eqref{e:strong}. In the right panel of Figure \ref{f:eig_plot} we plot 
$\beta$ and our estimator $\hat\beta$ as in \eqref{e:defhbb}. The estimator $\hat\beta$ provides a good approximation of $\beta$ provided that $\la>1$. In Figure \eqref{f:ev_plot} we plot 
$|\hat \bmv, \bmv\rangle|$, where the estimator $\hat\bmv$ is given as the left singular vector corresponding to the largest singular value of the unfolded matrix. Our theoretical prediction (blue curve) as in \eqref{e:strong1c} matches well with the the simulation for $\la>1$. For $\la\rightarrow 0$, our estimator behaves as poorly as random guess, i.e. taking $\hat\bmv$ as a random Gaussian vector (Green curve). For $\la$ in a small neighborhood of $1$, we don't have a good estimation of $|\langle \hat\bmv, \bmv\rangle|$, but only an upper bound \eqref{e:weak1c}.
 In the second panel of Figure \ref{f:ev_plot}, we zoom in around $\la=1$, the red curve, $\min\{n^{-1/6}, 1/\sqrt n |\la-1|\}$ corresponding to the  bound \eqref{e:weak1c}, provides a good upper bound of $|\langle \hat \bmv, \bmv\rangle|$.

\section{Low rank Perturbations of Long Random Matrices}\label{s:mainlong}

In this section we prove our main results as stated in Section \ref{s:main}. The proof of Theorem \ref{t:extreme} is given in Section \ref{s:longM}. We also collect the isotropic local law for long random matrices from \cite{alex2014isotropic}. 
In section \ref{s:master}, we derive a master equation which characterizes the outliers of the perturbed matrix $\beta \bmv \bmu^\top+ Z$. The master equation has been used intensively to study the low rank perturbation of random matrices, for both singular values and eigenvalues, see \cite{benaych2011eigenvalues, benaych2012singular, benaych2011fluctuations,huang2018mesoscopic }.
The proofs of Theorems \ref{t:eig} and \ref{t:ev} are given in Sections \ref{s:proveeig} and \ref{s:proveev} respectively.  

\subsection{Long Random Matrices}\label{s:longM}
Let $Z$ be an $n\times m$ random matrix, with entries satisfying Assumption \ref{a:Zasup}.  Let $\phi=\sqrt{m/n}$, with $n^{-C}\leq \phi\leq n^C$. In this section we recall some results of the sample covariance matrices in this setting from \cite{alex2014isotropic}. It turns out in this setting the correct normalization is to study $ZZ^\top/\phi$, which corresponds to the standard sample covariance matrices, with variance $\phi$. 

We denote the following $n$-dependent  Marchenko-Pastur law corresponding to the ratio $m/n=\phi^2\geq 1$, 
\begin{align}\label{e:MP}
\rho^{\phi}_{MP}(x)\rd x=\frac{\sqrt{\left(x-\left(\sqrt \phi-\frac{1}{\sqrt\phi}\right)^2\right)\left(\left(\sqrt \phi+\frac{1}{\sqrt\phi}\right)^2-x\right)}}{2\pi x/\phi}\rd x,
\end{align}
and its $1/n$-quanties as 
\begin{align}\label{e:quantiles}
\frac{i-1/2}{n}=\int_{\nu^\phi_i}^\infty \rho^\phi_{MP}(x)\rd x.
\end{align}
The normalization in \eqref{e:MP} is different from that in \eqref{e:MPa}, which corresponds to the sample covariance matrix $ZZ^\top$.
We remark that both the Marchenko-Pastur law $\phi_{MP}^\phi$ and its quantiles $\nu_i^\phi$ depend on $m,n$ through $\phi$.
We recall the following eigenvalue rigidity result from \cite[Theorem 2.10]{alex2014isotropic}.
\begin{theorem}[Eigenvalue Rigidity]\label{t:rigidity}
Under Assumption \ref{a:Zasup}, let $\phi=\sqrt{m/n}$ with  $n^{1/C}\leq \phi\leq n^C$, the eigenvalues $\la_1\geq\la_2\geq \cdots \geq \la_n\geq 0$ of $Z Z^\top/\phi$ are close to the quantiles of the Marchenko-Pastur law \eqref{e:MP}: fix any $\varepsilon>0$ and arbitrarily small $\delta>0$, with high probability it holds
\begin{align*}
|\la_i-\nu^\phi_i|\leq \frac{n^\delta}{i^{1/3}n^{2/3}},
\end{align*}
uniformly for all $i\in \qq{1, \floor{(1-\varepsilon)n}}$. If in addition $\phi-1\geq c$ for some constant $c>0$, then with high probability, we also have
\begin{align*}
|\la_i-\nu^\phi_i|\leq \frac{n^\delta}{(n+1-i)^{1/3}n^{2/3}},
\end{align*}
uniformly for all $i\in \qq{\floor{n/2},n}$, provided $n$ is large enough.
\end{theorem}

We denote the singular values of $Z$ as $s_1\geq s_2\geq \cdots \geq s_n\geq 0$. We can then restate Theorem \ref{t:rigidity} in terms of the singular values of $Z$, thanks to the following easy relation
\begin{align}\label{e:singeig}
\la_i=s_i^2/\phi,\quad 1\leq i\leq n.
\end{align}
Theorem \ref{t:extreme} is an easy consequence of Theorem \ref{t:rigidity} and the relation \eqref{e:singeig}.
\begin{proof}[Proof of Theorem \ref{t:extreme}]
The largest singular value $s_1$ of $Z$, and the largest eigenvalue $\la_1$ of $ZZ^\top/\phi$ are related by $\la_1=s_1^2/\phi$. Therefore, Theorem \ref{t:rigidity} implies that $|s_1^2/\phi-(\sqrt\phi+1/\sqrt \phi)^2|\leq n^\delta/n^{2/3}$, with high probability, provided $n$ is large enough. The claim \eqref{e:s1est} follows from rearranging.
\end{proof}

The empirical eigenvalue distribution of $ZZ^\top/\phi$ is close to the Marchenko-Pastur law \eqref{e:MP}. Thanks to the relation \eqref{e:singeig}, the empirical singular value distribution $\sum_{i=1}^n \delta_{s_i}/n$ of $Z$ is close to the push forward of the Marchenko-Pastur law \eqref{e:MP} (after proper normalization) by the map $x\mapsto \sqrt{x}$, 
\begin{align*}
&\phantom{{}={}}\rho_{\phi}(x):=2(x/\phi)\rho^\phi_{MP}(x^2/\phi)
=\frac{\sqrt{\left(x^2-\left(\phi-1\right)^2\right)\left(\left( \phi+1\right)^2-x^2\right)}}{\pi x}.
\end{align*}
We remark that $\rho_\phi$ is supported on $[\phi-1, \phi+1]$,

For later use, we denote the hermitization $\tilde Z$ of $Z$ as,
\begin{align}\label{e:tildeZ}
\tilde Z=\left[
\begin{array}{cc}
0 & Z\\
Z^\top & 0
\end{array}
\right].
\end{align}
Then $\tilde Z$ has $m-n$ zero eigenvalues, and its other eigenvalues are given by $\pm s_1, \pm s_2,\cdots, \pm s_n$. 
We denote the normalized Stieltjes transform of nonzero eigenvalues of $\tilde Z$ as
\begin{align*}
s(z)=\frac{1}{2n}\sum_{i=1}^n\frac{1}{z-s_i}+\frac{1}{z+s_i},
\end{align*}
and the Green's function of $\tilde Z$ as
\begin{align}\label{e:defG}
G(z)=\left(z
-\left[
\begin{array}{cc}
0 & Z\\
Z^\top & 0
\end{array}
\right]\right)^{-1}
=\left[
\begin{array}{cc}
(z-Z Z^\top/z)^{-1} & (z-Z Z^\top/z)^{-1}(Z/z)\\
(Z^\top/z)(z-Z Z^\top/z)^{-1} & (z-Z^\top Z/z)^{-1}
\end{array}
\right].
\end{align}

Then thanks to Theorem \ref{t:rigidity}, $s(z)$ is close to the Stieltjes transform of the symmetrized version of $\rho_\phi$,
\begin{align*}
m_\phi(z)=\frac{1}{2}\int\frac{\rho_\phi(x)+\rho_\phi(-x)}{z-x}\rd x
=\int\frac{z\rho_\phi(x)}{z^2-x^2}\rd x
=\int \frac{z\rho_{MP}^\phi(x)}{z^2-x\phi}\rd x.
\end{align*}
More explicitly, using the formula of $\rho_{MP}^\phi$ from \eqref{e:MP}, $m_\phi(z)$ is given by
\begin{align}\label{e:mpsi}
m_{\phi}(z)=\frac{-(\phi-1/\phi)+z^2/\phi\pm \sqrt{(z^2/\phi-(\sqrt\phi+1/\sqrt\phi)^2)(z^2/\phi-(\sqrt\phi-1/\sqrt\phi)^2)}}{2z/\phi},
\end{align}
and it satisfies the algebraic equation
\begin{align}\label{e:mgamma}
m_\phi^2+\left(\frac{\phi^2-1}{z}-z\right)m_\phi+1=0.
\end{align}


Let $z=E+\ri\eta$ with $\kappa\deq \min\{|E-(\phi-1)|, |E-(\phi+1)|\}$.
We denote the spectral domains
\begin{align}\begin{split}\label{e:defS}
&\bmS\deq \{z=E+\ri\eta: \kappa\leq \fc^{-1}, n^{-1+\fc}\leq \eta\leq \fc^{-1}\},\\
&\tilde\bmS \deq \{z=E+\ri\eta: E\not\in[(\phi-1),(\phi+1)], n^{-2/3+\fc}\leq \kappa\leq \fc^{-1}, 0<\eta\leq \fc^{-1}\}.
\end{split}\end{align}
The spectral domain $\bmS$ contains the spectral information inside the bulk, and the spectral domain $\tilde \bmS$ contains the spectral information close to the spectral edge.
We recall the following isotropic local law from \cite[Theorem 3.11, 3.12]{alex2014isotropic}, which will be used in Sections \ref{s:proveeig} and \ref{s:proveev}.
\begin{theorem}[Isotropic Local Law]\label{t:isotropiclaw}
For unit vectors $\bma\in \bR^n$ and $\bmb\in \bR^m$, with high probability, uniformly for $z\in \bmS$, we have
\begin{align*}\begin{split}
&\left| \langle \bma, (z- Z  Z^\top/z)^{-1}\bma\rangle-m_\phi(z)\langle \bma, \bma\rangle\right|\prec \sqrt{\frac{|\Im[m_\phi(z)]|}{n\eta}}+\frac{1}{n\eta},
\\
&\left| \langle \bma, (z- Z  Z^\top/z)^{-1}(Z/z)\bmb\rangle\right|\prec \frac{1}{\phi}\left(\sqrt{\frac{|\Im[m_\phi(z)]|}{n\eta}}+\frac{1}{n\eta}\right),
\\
&\left| \langle \bmb, (z- Z^\top  Z/z)^{-1}\bmb\rangle-\frac{\langle \bmb, \bmb\rangle}{z-m_\phi(z)}\right|\prec \frac{1}{\phi^2}\left(\sqrt{\frac{|\Im[m_\phi(z)]|}{n\eta}}+\frac{1}{n\eta}\right).
\end{split}
\end{align*}
Uniformly for $z\in \tilde \bmS$ we have the improved estimates:
\begin{align*}\begin{split}
&\left| \langle \bma, (z- Z  Z^\top/z)^{-1}\bma\rangle-m_\phi(z)\langle \bma, \bma\rangle\right|\prec \sqrt{\frac{|\Im[m_\phi(z)]|}{n\eta}},
\\
&\left| \langle \bma, (z- Z  Z^\top/z)^{-1}(Z/z)\bmb\rangle\right|\prec \frac{1}{\phi}\sqrt{\frac{|\Im[m_\phi(z)]|}{n\eta}},
\\
&\left| \langle \bmb, (z- Z^\top  Z/z)^{-1}\bmb\rangle-\frac{\langle \bmb, \bmb\rangle}{z-m_\phi(z)}\right|\prec \frac{1}{\phi^2}\sqrt{\frac{|\Im[m_\phi(z)]|}{n\eta}}.
\end{split}\end{align*}
\end{theorem}

\subsection{Master Equation}\label{s:master}
In this section, we derive a master equation, which characterize the outliers of the perturbed matrix $\beta \bmv \bmu^\top+ Z$. We denote the Hermitization of $\beta \bmv \bmu^\top+ Z$ as 
\begin{align}\label{e:hermit}
\left[
\begin{array}{cc}
\bmv/\sqrt 2 & \bmv\sqrt{2}\\
\bmu/\sqrt 2 & -\bmu\sqrt{2}
\end{array}
\right]
\left[
\begin{array}{cc}
\beta & 0\\
0  & -\beta
\end{array}
\right]
\left[
\begin{array}{cc}
\bmv/\sqrt 2 & \bmv\sqrt{2}\\
\bmu/\sqrt 2 & -\bmu\sqrt{2}
\end{array}
\right]^\top
+\left[
\begin{array}{cc}
0 & Z\\
Z^\top & 0
\end{array}
\right],
\end{align}
which encodes the spectral information of $\beta \bmv \bmu^\top+ Z$.
We can view \eqref{e:hermit} as a rank two perturbation of $\tilde Z$. We have the following well-known low rank perturbation formula 
\begin{lemma}\label{l:ABBA}
For two $(n+m)\times r$ matrices $A$ and $B$, it holds that
\begin{align}\label{e:sdt}
\det(I-AB^\top)=\det(I-B^\top A).
\end{align}
The matrix $I-AB^\top$ is invertible if and only if $I-B^\top A$ is invertible, and 
\begin{align}\label{e:minv}
(I-AB^\top)^{-1}=I+A(I-B^\top A)^{-1} B^\top.
\end{align}
\end{lemma}
\begin{proof}
The identity \eqref{e:sdt} is Sylvester's determinant theorem. And the second identity \eqref{e:minv} is known as the matrix inversion lemma or Woodbury's matrix identity. 
\end{proof}

We will use Lemma \ref{l:ABBA} to study \eqref{e:hermit}, which is a rank two perturbation of $\tilde Z$. Its eigenvalues are given by the roots of the characteristic polynomials
\begin{align*}
&\phantom{{}={}}\det\left(z-\left[
\begin{array}{cc}
\bmv/\sqrt 2 & \bmv\sqrt{2}\\
\bmu/\sqrt 2 & -\bmu\sqrt{2}
\end{array}
\right]
\left[
\begin{array}{cc}
\beta & 0\\
0  & -\beta
\end{array}
\right]
\left[
\begin{array}{cc}
\bmv/\sqrt 2 & \bmv\sqrt{2}\\
\bmu/\sqrt 2 & -\bmu\sqrt{2}
\end{array}
\right]^\top
-\left[
\begin{array}{cc}
0 & Z\\
Z^\top & 0
\end{array}
\right]\right)\\
&=\det\left(z
-\left[
\begin{array}{cc}
0 & Z\\
Z^\top & 0
\end{array}
\right]\right)\det\left(I-G(z)\left[
\begin{array}{cc}
\bmv/\sqrt 2 & \bmv\sqrt{2}\\
\bmu/\sqrt 2 & -\bmu\sqrt{2}
\end{array}
\right]
\left[
\begin{array}{cc}
\beta & 0\\
0  & -\beta
\end{array}
\right]
\left[
\begin{array}{cc}
\bmv/\sqrt 2 & \bmv\sqrt{2}\\
\bmu/\sqrt 2 & -\bmu\sqrt{2}
\end{array}
\right]^\top
\right)\\
&=\det\left(z
-\left[
\begin{array}{cc}
0 & Z\\
Z^\top & 0
\end{array}
\right]\right)\det\left(I-U^\top A(z) 
\right),
\end{align*}
where we used Lemma \ref{l:ABBA} for $r=2$, and
\begin{align}\label{e:AU}
A(z):=G(z)\left[
\begin{array}{cc}
\bmv/\sqrt 2 & \bmv\sqrt{2}\\
\bmu/\sqrt 2 & -\bmu\sqrt{2}
\end{array}
\right]
\left[
\begin{array}{cc}
\beta & 0\\
0  & -\beta
\end{array}
\right],\quad
U:=
\left[
\begin{array}{cc}
\bmv/\sqrt 2 & \bmv\sqrt{2}\\
\bmu/\sqrt 2 & -\bmu\sqrt{2}
\end{array}
\right].
\end{align}
Therefore, the singular values of $\beta \bmu\bmv^\top+Z$ (which are not eigenvalues of $\tilde Z$) are characterized by 
\begin{align}\label{e:IUA}
\det(I-U^\top A(z))=0.
\end{align}
The equation \eqref{e:IUA} can be used to characterize the outliers of $\beta \bmv\bmu^\top+Z$, and we will use it to prove Theorem \ref{t:eig} in Section \ref{s:proveeig}.

Thanks to Lemma \ref{l:ABBA}, for $z\in \bC^+$ in the upper half plane, we can write explicitly the Green's of the Hermitization \eqref{e:hermit} of $\beta \bmv\bmu^\top +Z$,
\begin{align}\begin{split}\label{e:Greena}
&\phantom{{}={}}\left(z-\left[
\begin{array}{cc}
\bmv/\sqrt 2 & \bmv\sqrt{2}\\
\bmu/\sqrt 2 & -\bmu\sqrt{2}
\end{array}
\right]
\left[
\begin{array}{cc}
\beta & 0\\
0  & -\beta
\end{array}
\right]
\left[
\begin{array}{cc}
\bmv/\sqrt 2 & \bmv\sqrt{2}\\
\bmu/\sqrt 2 & -\bmu\sqrt{2}
\end{array}
\right]^\top
-\left[
\begin{array}{cc}
0 & Z\\
Z^\top & 0
\end{array}
\right]\right)^{-1}\\
&=
\left(I-G(z)\left[
\begin{array}{cc}
\bmv/\sqrt 2 & \bmv\sqrt{2}\\
\bmu/\sqrt 2 & -\bmu\sqrt{2}
\end{array}
\right]
\left[
\begin{array}{cc}
\beta & 0\\
0  & -\beta
\end{array}
\right]
\left[
\begin{array}{cc}
\bmv/\sqrt 2 & \bmv\sqrt{2}\\
\bmu/\sqrt 2 & -\bmu\sqrt{2}
\end{array}
\right]^\top
\right)^{-1}G(z)\\
&=\left(I+A(z)(I-U^\top A(z))^{-1}U^\top\right)G(z).
\end{split}\end{align}
The Green's function \eqref{e:Greena} contains the information of eigenvectors, and can be used to study the singular vectors of the outliers of $\beta \bmv\bmu^\top+Z$.  And we will use it to prove Theorem \ref{t:ev} in Section \ref{s:proveev}.

\subsection{Proof of Theorem \ref{t:eig}}\label{s:proveeig}
Let $\phi=\sqrt{m/n}$. We denote the singular values of $\beta \bmv\bmu^\top+Z$ as $\hat s_1\geq \hat s_2\geq \cdots\geq  \hat s_n\geq 0$, with corresponding normalized left and right singular vectors as $\hat \bmv_1, \hat\bmv_2, \cdots, \hat \bmv_n\in \bR^n$, and $\hat \bmu_1, \hat\bmu_2, \cdots, \hat \bmu_n\in \bR^m$. Then the nonzero eigenvalues of its Hermitization 
\begin{align}\label{e:hermitc}
\left[
\begin{array}{cc}
\bmv/\sqrt 2 & \bmv\sqrt{2}\\
\bmu/\sqrt 2 & -\bmu\sqrt{2}
\end{array}
\right]
\left[
\begin{array}{cc}
\beta & 0\\
0  & -\beta
\end{array}
\right]
\left[
\begin{array}{cc}
\bmv/\sqrt 2 & \bmv\sqrt{2}\\
\bmu/\sqrt 2 & -\bmu\sqrt{2}
\end{array}
\right]^\top
+\left[
\begin{array}{cc}
0 & Z\\
Z^\top & 0
\end{array}
\right],
\end{align}
are given by $\hat s_1, -\hat s_1, \hat s_2, -\hat s_2,\cdots, \hat s_n, -\hat s_n$. The corresponding normalized eigenvectors are given by
$[\hat\bmv_1, \hat\bmu_1]/\sqrt 2, [\hat\bmv_1, -\hat\bmu_1]/\sqrt 2, [\hat\bmv_2, \hat\bmu_2]/\sqrt 2, [\hat\bmv_2, -\hat\bmu_2]/\sqrt 2,\cdots, [\hat\bmv_n, \hat\bmu_n]/\sqrt 2, [\hat\bmv_n, -\hat\bmu_n]/\sqrt 2$.

We recall from Theorem \ref{t:rigidity} that with high probability the singular values of $Z$ are bounded by $(\phi+1)+n^{-2/3+\fc/2}$, i.e. $s_1\leq  (\phi+1)+n^{-2/3+\fc/2}$. In the remaining of this section, we restrict ourselves to the event that $s_1\leq  (\phi+1)+n^{-2/3+\fc/2}$, which holds with high probability. We can view \eqref{e:hermitc} as a rank two perturbation of $\tilde Z$,
\begin{align}\label{e:rankone}
\beta \left[
\begin{array}{c}
\bmv/\sqrt 2 \\
\bmu/\sqrt 2 
\end{array}
\right]
\left[
\begin{array}{c}
\bmv/\sqrt 2 \\
\bmu/\sqrt 2
\end{array}
\right]^\top
-\beta \left[
\begin{array}{c}
\bmv/\sqrt 2 \\
-\bmu/\sqrt 2 
\end{array}
\right]
\left[
\begin{array}{c}
\bmv/\sqrt 2 \\
-\bmu/\sqrt 2
\end{array}
\right]^\top
+\left[
\begin{array}{cc}
0 & Z\\
Z^\top & 0
\end{array}
\right],
\end{align}
By the variational formula of eigenvalues, $\hat s_2$ is upper bounded by the second largest eigenvalue of 
\begin{align}\label{e:rankone}
\beta \left[
\begin{array}{c}
\bmv/\sqrt 2 \\
\bmu/\sqrt 2 
\end{array}
\right]
\left[
\begin{array}{c}
\bmv/\sqrt 2 \\
\bmu/\sqrt 2
\end{array}
\right]^\top
+\left[
\begin{array}{cc}
0 & Z\\
Z^\top & 0
\end{array}
\right].
\end{align}
Again the second largest eigenvalue of \eqref{e:rankone} is upper bounded by $s_1$. Thus we have $\hat s_2\leq s_1\leq (\phi+1)+n^{-2/3+\fc/2}$. Therefore, \eqref{e:hermitc} can have at most one outlier eigenvalue.

We recall from \eqref{e:IUA}, the matrix $\beta \bmv\bmu^\top+Z$ has an outlier singular value  bigger than $ (\phi+1)+n^{-2/3+\fc}$, if and only if 
\begin{align}\label{e:IUAc}
\det(I-U^\top A(x))=0
\end{align}
has an zero with $x\geq  (\phi+1)+n^{-2/3+\fc}$, where
\begin{align*}
A(z)=G(z)\left[
\begin{array}{cc}
\bmv/\sqrt 2 & \bmv\sqrt{2}\\
\bmu/\sqrt 2 & -\bmu\sqrt{2}
\end{array}
\right]
\left[
\begin{array}{cc}
\beta & 0\\
0  & -\beta
\end{array}
\right],\quad
U=
\left[
\begin{array}{cc}
\bmv/\sqrt 2 & \bmv\sqrt{2}\\
\bmu/\sqrt 2 & -\bmu\sqrt{2}
\end{array}
\right].
\end{align*}

We can rewrite the equation \eqref{e:IUAc} as
\begin{align}\label{e:eigeq}
&\det\left(\left[
\begin{array}{cc}
1/\beta & 0\\
0  & -1/\beta
\end{array}
\right]-\left[
\begin{array}{cc}
\bmv/\sqrt 2 & \bmv\sqrt{2}\\
\bmu/\sqrt 2 & -\bmu\sqrt{2}
\end{array}
\right]^\top G(x)\left[
\begin{array}{cc}
\bmv/\sqrt 2 & \bmv\sqrt{2}\\
\bmu/\sqrt 2 & -\bmu\sqrt{2}
\end{array}
\right]
\right)=0.
\end{align}


We recall the Green's function $G(x)$ from \eqref{e:defG}
\begin{align*}
G(x)=\left(x
-\left[
\begin{array}{cc}
0 & Z\\
Z^\top & 0
\end{array}
\right]\right)^{-1}
=\left[
\begin{array}{cc}
(x-Z Z^\top/x)^{-1} & (x-Z Z^\top/x)^{-1}(Z/x)\\
(Z^\top/x)(x-Z Z^\top/x)^{-1} & (x-Z^\top Z/x)^{-1}
\end{array}
\right],
\end{align*}
and denote the quantities
\begin{align}\begin{split}\label{e:defABC}
&\cA(x)= \langle \bmv, (x- Z  Z^\top/x)^{-1}\bmv\rangle,\\
&\cB(x)=\langle \bmv, (x- Z  Z^\top/x)^{-1}(Z/x)\bmu\rangle,\\
&\cC(x)=\langle \bmu, (x- Z^\top  Z/x)^{-1}\bmu\rangle
\end{split}\end{align}
With these notations, we can rewrite \eqref{e:eigeq} as
\begin{align}\label{e:repll}
\det\left(\left[
\begin{array}{cc}
1/\beta & 0\\
0  & -1/\beta
\end{array}
\right]-\frac{1}{2}\left[
\begin{array}{cc}
\cA(x)+2\cB(x)+\cC(x) & \cA(x)-\cC(x)\\
\cA(x)-\cC(x) & \cA(x)-2\cB(x)+\cC(x)
\end{array}
\right]
\right)=0.
\end{align}
It simplifies to 
\begin{align}\label{e:eq0}
\left(\frac{1}{\beta}-\cB(x)\right)^2=\cA(x)\cC(x).
\end{align}

Let $x-(\phi+1)=\kappa\geq n^{-2/3+\fc}$, then $x\in \tilde \bmS$ as defined in \eqref{e:defS}.
Thanks to the Square root behavior \eqref{e:mpsi} of $m_\phi(z)$ around the spectral edge $z=\phi+1$, we have
\begin{align*}
|\Im[m_\phi(z)]|\asymp \frac{\eta}{\sqrt{\kappa+\eta}}, \quad z=\phi+1+\kappa+\ri\eta,\quad \kappa\geq 0.
\end{align*}
In particular we have 
\begin{align}\label{e:m/eta}
\sqrt{\frac{|\Im[m_{\phi}(z)]|}{n\eta}}\asymp \frac{1}{n^{1/2}\kappa^{1/4}},\quad z=\phi+1+\kappa+\ri\eta,\quad \kappa\geq 0.
\end{align}
Thanks to Theorem \ref{t:isotropiclaw} with $\bma=\bmv$ and $\bmb=\bmu$, by plugging \eqref{e:m/eta}, we have
\begin{align}\begin{split}\label{e:ABCcopyhi}
&\cA(x)=m_\phi(x)+\OO_\prec\left(\frac{1}{n^{1/2}\kappa^{1/4}}\right),\\
&\cB(x)=\OO_\prec\left(\frac{1}{\phi n^{1/2}\kappa^{1/4}}\right),\\
&\cC(x)=\frac{1}{x-m_\phi(x)}+\OO_\prec\left(\frac{1}{\phi^2 n^{1/2}\kappa^{1/4}}\right),
\end{split}\end{align}
where the error terms are continuous in $x$. 
Then by plugging \eqref{e:ABCcopyhi} we can rewrite \eqref{e:eq0} as
\begin{align}\label{e:eq1}
\frac{1}{\beta^2}= \frac{m_\phi(x)}{x-m_\phi(x)} +\OO_\prec\left(\frac{1}{\phi n^{1/2}\kappa^{1/4}}\right).
\end{align}
We recall the algebraic equation of $m_\phi$ from \eqref{e:mgamma}
\begin{align}\label{e:meq}
m_\phi^2+\left(\frac{\phi^2-1}{x}-x\right)m_\phi+1=0,
\end{align}

By further rearranging \eqref{e:eq1} we get
\begin{align}\label{e:eqha}
\frac{m_\phi(x)}{x-m_\phi(x)}=\frac{1}{\beta^2}+\OO_\prec\left(\frac{1}{\phi n^{1/2}\kappa^{1/4}}\right).
\end{align}
Since $m_\phi(x)$ is monotone decreasing for $x\geq \phi+1$, the lefthand side is monotone decreasing for $x\geq \phi+1$. Using the formula \eqref{e:mpsi} for $m_\phi(x)$, for $x=\phi+1$, $m_\phi(x)=1$. For $x\rightarrow \infty$, we have $m_\phi(x)\rightarrow 0$.  And for $x=\phi+1+\kappa$, with $\kappa=\oo(1)$, we have 
\begin{align*}
m_\phi(x)=1-C\sqrt \kappa+\OO(\kappa),
\end{align*}
for some constant $C>0$. In this regime, the lefthand side of \eqref{e:eqha} behaves 
\begin{align}\label{e:mphidiv}
\frac{m_\phi(x)}{x-m_\phi(x)}=\frac{1-C\sqrt \kappa +\OO(\kappa)}{\phi+C\sqrt \kappa +\OO(\kappa)}
=\frac{1}{\phi}\left(1-(C+C/\phi)\sqrt\kappa+\OO(\kappa)\right)
\end{align}

We take 
\begin{align*}
\beta=\la\sqrt \phi.
\end{align*}
In the following we first prove \eqref{e:weak}.
\begin{proof}[Proof of \eqref{e:weak}]
If $\la\leq 1+n^{-1/3+\fc}$, in the following we show that \eqref{e:eq1} does not have solution with $x\geq \phi+1+n^{-2/3+3\fc}$. For $x\geq \phi+1+n^{-2/3+3\fc}$, using \eqref{e:mphidiv} and the monotonicity of $m_\phi(x)$, we have
\begin{align}\begin{split}\label{e:left}
\frac{m_\phi(x)}{x-m_\phi(x)}
&\leq \frac{m_\phi(\phi+1+n^{-2/3+3\fc})}{\phi+1+n^{-2/3+3\fc}-m_\phi(\phi+1+n^{-2/3+3\fc})}\\
&=\frac{1}{\phi}(1-(C+C/\phi)n^{-1/3+3\fc/2}+\OO(n^{-2/3+3\fc})).
\end{split}\end{align}
For $x\geq \phi+1+n^{-2/3+3\fc}$, the righthand side of \eqref{e:eqha} is lower bounded by
\begin{align}\begin{split}\label{e:right}
&\phantom{{}={}}\frac{1}{\beta^2}+\OO_\prec\left(\frac{1}{\phi n^{1/2}\kappa^{1/4}}\right)
= \frac{1}{\la^2\phi}+\OO_\prec\left(\frac{1}{\phi n^{1/2}\kappa^{1/4}}\right)\\
&\geq \frac{1}{\phi}(1-2n^{-1/3+\fc}+\OO(n^{-2/3+2\fc}))+\OO_\prec\left( \frac{1}{\phi n^{1/3+3\fc/4}}\right)
>\frac{m_\phi(x)}{x-m_\phi(x)}.
\end{split}\end{align}
Therefore, \eqref{e:eqha} does not have solution for $x\geq \phi+1+n^{-2/3+3\fc}$. We conclude that $\hat s_1\leq \phi+1+n^{-2/3+3\fc}$. This finishes the proof \eqref{e:weak}.
\end{proof}

In the following, we study the case that $\la\geq 1+n^{-1/3+\fc}$. 
\begin{proof}[Proof of \eqref{e:strong}]
Similarly to \eqref{e:left} and \eqref{e:right} we have for $x= \phi+1+n^{-2/3+\fc}$
\begin{align*}
\frac{1}{\beta^2}+\OO_\prec\left(\frac{1}{\phi n^{1/2}\kappa^{1/4}}\right)\leq \frac{m_\phi(x)}{x-m_\phi(x)}.
\end{align*}
And for $x\rightarrow \infty$,
\begin{align*}
\lim_{x\rightarrow \infty }\frac{1}{\beta^2}+\OO_\prec\left(\frac{1}{\phi n^{1/2}\kappa^{1/4}}\right)=\frac{1}{\beta^2}\geq  \lim_{x\rightarrow \infty }\frac{m_\phi(x)}{x-m_\phi(x)}=0.
\end{align*}
Therefore \eqref{e:eqha} a has a solution in the interval $[\phi+1+n^{-2/3+\fc}, \infty)$, which is $\hat s_1$. 
In particularly we have $\hat s_1\geq \phi+1+n^{-2/3+\fc}$, and $\beta \bmv\bmu^\top +Z$ has an outlier. In the following, we compute the value of $\hat s_1$.

We rewrite the equation \eqref{e:eqha} as
\begin{align}\label{e:mphi}
\frac{m_\phi(x)}{x-m_\phi(x)}= \frac{\tau}{\phi},\quad
\tau=\frac{1}{\la^2}+\OO_\prec\left(\frac{1}{n^{1/2}\kappa^{1/4}}\right).
\end{align}

We can solve for $x$ and $m_\phi(x)$ using \eqref{e:meq} and \eqref{e:mphi}
\begin{align}\begin{split}\label{e:xmexp}
&x=\sqrt{\phi^2+(\tau+1/\tau)\phi+1}
=\sqrt{\phi^2+(\la^2+1/\la^2)\phi+1}+\OO_\prec \left(\frac{(\la-1)}{\phi}\frac{1}{n^{1/2}\kappa^{1/4}}+\frac{1}{\phi}\frac{1}{n\kappa^{1/2}}\right),\\
&m_\phi=\frac{x}{1+\phi/\tau}=\left(1+\OO_\prec\left(\frac{1}{n^{1/2}\kappa^{1/4}}\right)\right) \sqrt{\frac{\phi+\la^2}{\la^2(1+\la^2\phi)}}.
\end{split}\end{align}
We recall that
\begin{align*}
x=\phi+1+\kappa,\quad \kappa\geq n^{-2/3+\fc}
\end{align*}
By a Taylor expansion of the first estimate in \eqref{e:xmexp}, for $\la=1+\oo(1)$, there exists some constant $C>0$,
\begin{align*}
x\geq \phi+1+C(\la-1)^2+\OO_\prec \left(\frac{(\la-1)}{\phi}\frac{1}{n^{1/2}\kappa^{1/4}}+\frac{1}{\phi}\frac{1}{n\kappa^{1/2}}\right)
\end{align*}
By our assumption that $\la\geq 1+n^{-1/3+\fc}$. Then we have $\kappa=x-\phi-1\asymp (\la-1)^2\gtrsim n^{-2/3+2\fc}$. We can use this estimate of $\kappa$ to simplify the error terms in \eqref{e:xmexp}.
In summary we have for $\beta=\la\sqrt \phi$ and $\la\geq 1+n^{-1/3+\fc}$, $\beta \bmv\bmu^\top + Z$ has an outlier singular value at
\begin{align*}
x=\sqrt{\phi^2+(\la^2+1/\la^2)\phi+1}+\OO_\prec \left(\frac{(\la-1)^{1/2}}{n^{1/2} \phi}\right)
\end{align*}
This finishes the proof of \eqref{e:strong}. 
\end{proof}


\subsection{Proof of Theorem \ref{t:ev}}\label{s:proveev}


Let $\phi=\sqrt{m/n}$. We first study the case that $\la\geq 1+n^{-1/3+\fc}$. 

\begin{proof}[Proof of \eqref{e:strong1} and \eqref{e:strong22}]
If $\la\geq 1+n^{-1/3+\fc}$, Theorem \eqref{t:eig} implies that $\beta \bmv\bmu^\top +Z$ has an outlier singular value $x=\hat s_1=\phi+1+\kappa$ with $\kappa\asymp (\la-1)^2$. Therefore, $x$ is an eigenvalue of the Hermitized matrix:
\begin{align}\label{e:eigvvw}
x \bmw=\left(\left[
\begin{array}{cc}
\bmv/\sqrt 2 & \bmv\sqrt{2}\\
\bmu/\sqrt 2 & -\bmu\sqrt{2}
\end{array}
\right]
\left[
\begin{array}{cc}
\beta & 0\\
0  & -\beta
\end{array}
\right]
\left[
\begin{array}{cc}
\bmv/\sqrt 2 & \bmv\sqrt{2}\\
\bmu/\sqrt 2 & -\bmu\sqrt{2}
\end{array}
\right]^\top
+\left[
\begin{array}{cc}
0 & Z\\
Z^\top & 0
\end{array}
\right]\right)\bmw,
\end{align}
for some unit vector $\bmw=[\hat\bmv_1, \hat\bmu_1]^\top/\sqrt 2$, where $\hat\bmv_1, \hat\bmu_1$ are the left and right singular vector of $\beta \bmv\bmu^\top +Z$ corresponding to the singular value $\hat s_1$. 
By rearranging \eqref{e:eigvvw}, we get
\begin{align}\label{e:weq}
\bmw=\left(x-\left[
\begin{array}{cc}
0 & Z\\
Z^\top & 0
\end{array}
\right]\right)^{-1}\left[
\begin{array}{cc}
\bmv/\sqrt 2 & \bmv\sqrt{2}\\
\bmu/\sqrt 2 & -\bmu\sqrt{2}
\end{array}
\right]
\left[
\begin{array}{cc}
\beta & 0\\
0  & -\beta
\end{array}
\right]
\left[
\begin{array}{cc}
\bmv/\sqrt 2 & \bmv\sqrt{2}\\
\bmu/\sqrt 2 & -\bmu\sqrt{2}
\end{array}
\right]^\top\bmw.
\end{align}
We denote the length two vector which is the projection of $\bmw$ on $\bmv$ and $\bmu$ direction
\begin{align*}
\tilde \bmw=\left[
\begin{array}{cc}
\bmv/\sqrt 2 & \bmv\sqrt{2}\\
\bmu/\sqrt 2 & -\bmu\sqrt{2}
\end{array}
\right]^\top\bmw=\frac{1}{2}\left[
\begin{array}{cc}
\bmv & \bmv\\
\bmu & -\bmu
\end{array}
\right]^\top
\left[
\begin{array}{c}
\hat \bmv_1\\
\hat\bmu_1
\end{array}
\right],
\end{align*}
then it satisfies 
\begin{align}\label{e:ttddw}
\tilde \bmw
=\left[
\begin{array}{cc}
\bmv/\sqrt 2 & \bmv\sqrt{2}\\
\bmu/\sqrt 2 & -\bmu\sqrt{2}
\end{array}
\right]^\top\left(x-\left[
\begin{array}{cc}
0 & Z\\
Z^\top & 0
\end{array}
\right]\right)^{-1}\left[
\begin{array}{cc}
\bmv/\sqrt 2 & \bmv\sqrt{2}\\
\bmu/\sqrt 2 & -\bmu\sqrt{2}
\end{array}
\right]
\left[
\begin{array}{cc}
\beta & 0\\
0  & -\beta
\end{array}
\right]
\tilde\bmw.
\end{align}
We denote $\tilde \bmw=[\tilde \bmw_1, \tilde \bmw_2]^\top$ with $\tilde \bmw_1\in \bR^n$ and $\tilde \bmw_2\in \bR^m$.  
The inner products $\langle \hat{\bmv}_1, \bmv\rangle$ and $\langle \hat{\bmu}_1, \bmu\rangle$ are given by
\begin{align}\label{e:innerp}
\langle \hat{\bmv}_1, \bmv\rangle =\tilde \bmw_1+\tilde \bmw_2,\quad
\langle \hat{\bmu}_1, \bmu\rangle =\tilde \bmw_1-\tilde \bmw_2.
\end{align}

By taking norm on both sides of \eqref{e:weq}, we get another equation for $\tilde w$,
\begin{align}\label{e:wGw}
1=\tilde \bmw^\top
\left[
\begin{array}{cc}
\beta & 0\\
0  & -\beta
\end{array}
\right]\left[
\begin{array}{cc}
\bmv/\sqrt 2 & \bmv\sqrt{2}\\
\bmu/\sqrt 2 & -\bmu\sqrt{2}
\end{array}
\right]^\top
\left(x-\left[
\begin{array}{cc}
0 & Z\\
Z^\top & 0
\end{array}
\right]\right)^{-2}\left[
\begin{array}{cc}
\bmv/\sqrt 2 & \bmv\sqrt{2}\\
\bmu/\sqrt 2 & -\bmu\sqrt{2}
\end{array}
\right]
\left[
\begin{array}{cc}
\beta & 0\\
0  & -\beta
\end{array}
\right]
\tilde\bmw.
\end{align}

Using the notations from \eqref{e:defABC}, we can rewrite \eqref{e:ttddw} as 
\begin{align}\label{e:evrelat}
\tilde \bmw=\frac{1}{2}\left[
\begin{array}{cc}
\cA(x)+2\cB(x)+\cC(x) & \cA(x)-\cC(x)\\
\cA(x)-\cC(x) & \cA(x)-2\cB(x)+\cC(x)
\end{array}
\right]\left[
\begin{array}{cc}
\beta & 0\\
0  & -\beta
\end{array}
\right]
\tilde \bmw,
\end{align}
where we recall from \eqref{e:eq0}
\begin{align}\label{e:relate}
\left(\frac{1}{\beta}-\cB(x)\right)^2=\cA(x)\cC(x).
\end{align}
Thanks to Theorem \ref{t:isotropiclaw},  for $x=\phi+1+\kappa$ with $\kappa\geq n^{-2/3+\fc}$, we have
\begin{align}\begin{split}\label{e:ABCbb}
&\cA(x)=m_\phi(x)+\OO_\prec\left(\frac{1}{n^{1/2}\kappa^{1/4}}\right),\\
&\cB(x)=\OO_\prec\left(\frac{1}{\phi n^{1/2}\kappa^{1/4}}\right)),\\
&\cC(x)=\frac{1}{x-m_\phi(x)}+\OO_\prec\left(\frac{1}{\phi^2 n^{1/2}\kappa^{1/4}}\right),
\end{split}\end{align}
where we used \eqref{e:m/eta}.

On the event that the singular values of $Z$ are bounded by $(\phi+1)+n^{-2/3+\fc/2}$, i.e. $s_1\leq  (\phi+1)+n^{-2/3+\fc/2}$, we have that $\cA(z), \cB(z),\cC(z)$ are analytic for $\Re[z]\geq (\phi+1)+n^{-2/3+\fc/2}$. We take a  contour $\omega=\{z: |z-x|=\kappa/2\}$. Inside the contour $\Re[z]\geq \phi+1+n^{-2/3+\fc}/2$, $\cA(z)$ is analytic. 
Then we can rewrite $\cA'(x)$ as a contour integral
\begin{align}\begin{split}\label{e:ABC'0}
\cA'(x)&=\frac{1}{2\pi\ri}\oint_\omega\frac{\cA(z)}{(z-x)^2}\rd z
=\frac{1}{2\pi\ri}\oint_\omega\frac{m_\phi(z)}{(z-x)^2}\rd z
+\OO_\prec\left(\oint_\omega \frac{\rd z}{|z-x|^2}\frac{1}{n^{1/2}\kappa^{1/4}}\right)\\
&=m'_\phi(x)+\OO_\prec\left(\frac{1}{n^{1/2}\kappa^{5/4}}\right),
\end{split}\end{align}
where in the last equality, we used that the total length of the contour $\omega$ is of order $\kappa$. Similar argument also gives us the estimates of $\cB'(x)$ and $\cC'(x)$,
\begin{align}\begin{split}\label{e:ABC'}
&\cB'(x)=\OO_\prec\left(\frac{1}{\phi n^{1/2}\kappa^{5/4}}\right),\\
&\cC'(x)=\frac{m'_\phi(x)-1}{(x-m_\phi(x))^2}+\OO_\prec\left(\frac{1}{\phi^2 n^{1/2}\kappa^{5/4}}\right).
\end{split}\end{align}

By slightly rearranging \eqref{e:evrelat}, we get that 
\begin{align*}
&\phantom{{}={}}\left[
\begin{array}{cc}
\beta & 0\\
0  & -\beta
\end{array}
\right]\tilde \bmw
=\frac{1}{2}\left[
\begin{array}{cc}
\beta & 0\\
0  & -\beta
\end{array}
\right]\left[
\begin{array}{cc}
\cA(x)+2\cB(x)+\cC(x) & \cA(x)-\cC(x)\\
\cA(x)-\cC(x) & \cA(x)-2\cB(x)+\cC(x)
\end{array}
\right]\left[
\begin{array}{cc}
\beta & 0\\
0  & -\beta
\end{array}
\right]
\tilde \bmw.
\end{align*}
Therefore 
\begin{align*}
\left[
\begin{array}{cc}
\beta & 0\\
0  & -\beta
\end{array}
\right]
\tilde \bmw,
\end{align*}
is an eigenvector of the following matrix with eigenvalue $1$,
\begin{align}
\frac{1}{2}\left[
\begin{array}{cc}
\beta & 0\\
0  & -\beta
\end{array}
\right]\left[
\begin{array}{cc}
\cA(x)+2\cB(x)+\cC(x) & \cA(x)-\cC(x)\\
\cA(x)-\cC(x) & \cA(x)-2\cB(x)+\cC(x)
\end{array}
\right].\label{e:themm}
\end{align}

By plugging \eqref{e:relate} into \eqref{e:themm}, we can rewrite it as
\begin{align*}
&I+\frac{\beta}{2} \left[
\begin{array}{c}
\sqrt{\cA(x)}-\sqrt{\cC(x)}\\
-(\sqrt{\cA(x)}+\sqrt{\cC(x)})
\end{array}
\right]
 \left[
\begin{array}{cc}
\sqrt{\cA(x)}-\sqrt{\cC(x)}&
\sqrt{\cA(x)}+\sqrt{\cC(x)}
\end{array}
\right],
\end{align*}
which has an eigenvector $[\sqrt{\cC(x)}+\sqrt{\cA(x)}, \sqrt{\cC(x)}-\sqrt{\cA(x)}]^\top$ with eigenvalue $1$. We conclude that 
the eigenvector $\tilde \bmw$ satisfies
\begin{align}\label{e:twexxq}
\left[
\begin{array}{cc}
\beta & 0\\
0  & -\beta
\end{array}
\right]
\tilde\bmw=c 
\left[
\begin{array}{c}
\sqrt{\cC(x)}+\sqrt{\cA(x)}\\
\sqrt{\cC(x)}-\sqrt{\cA(x)}.
\end{array}
\right]
\end{align}

We need to use \eqref{e:wGw} to determine $c$ in the above expression, 
\begin{align*}
&\phantom{{}={}}\left[
\begin{array}{cc}
\bmv/\sqrt 2 & \bmv\sqrt{2}\\
\bmu/\sqrt 2 & -\bmu\sqrt{2}
\end{array}
\right]^\top\left(x-\left[
\begin{array}{cc}
0 & Z\\
Z^\top & 0
\end{array}
\right]\right)^{-2}\left[
\begin{array}{cc}
\bmv/\sqrt 2 & \bmv\sqrt{2}\\
\bmu/\sqrt 2 & -\bmu\sqrt{2}
\end{array}
\right]\\
&=-\del_x \left[
\begin{array}{cc}
\bmv/\sqrt 2 & \bmv\sqrt{2}\\
\bmu/\sqrt 2 & -\bmu\sqrt{2}
\end{array}
\right]^\top\left(x-\left[
\begin{array}{cc}
0 & Z\\
Z^\top & 0
\end{array}
\right]\right)^{-1}\left[
\begin{array}{cc}
\bmv/\sqrt 2 & \bmv\sqrt{2}\\
\bmu/\sqrt 2 & -\bmu\sqrt{2}
\end{array}
\right]\\
&=
-\frac{1}{2}\left[
\begin{array}{cc}
\cA'(x)+2\cB'(x)+\cC'(x) & \cA'(x)-\cC'(x)\\
\cA'(x)-\cC'(x) & \cA'(x)-2\cB'(x)+\cC'(x)
\end{array}
\right]
\end{align*}
By plugging the above expression to \eqref{e:wGw} and using \eqref{e:ABCbb}, \eqref{e:ABC'0}, \eqref{e:ABC'} to simplify, we get
\begin{align}\begin{split}\label{e:computec}
1&=-2c^2(\cC(x)\cA'(x)+\cA(x) \cC'(x)+2\sqrt{\cA(x)\cC(x)}\cB'(x))\\
&=-2c^2\left( \frac{m_\phi'(x)}{x-m_\phi(x)}+\frac{m_\phi(m'_\phi-1)}{(x-m_\phi)^2}+\OO_\prec\left(\frac{1}{\phi n^{1/2}\kappa^{5/4}}\right)\right)\\
&=-2c^2\left( \frac{m_\phi'(x)}{x-m_\phi(x)}+\frac{m_\phi(m'_\phi-1)}{(x-m_\phi)^2}\right)\left(1+\OO_\prec\left(\frac{1}{ n^{1/2}\kappa^{3/4}}\right)\right),
\end{split}\end{align}
where we also used that $|m_\phi'(x)|\asymp 1/\sqrt\kappa$, where $\kappa=x-(\phi+1)$.

We recall the inner products $\langle \hat{\bmv}_1, \bmv\rangle$ and $\langle \hat{\bmu}_1, \bmu\rangle$ from \eqref{e:innerp}, using \eqref{e:twexxq} and \eqref{e:ABCcopyhi}
\begin{align}\label{e:innerp2}
&\langle \hat{\bmv}_1, \bmv\rangle =\tilde \bmw_1+\tilde \bmw_2=
\frac{2c}{\beta}\sqrt{\cA(x)}
=\frac{2c \sqrt{m_\phi}}{\beta}\left(1+\OO_\prec\left(\frac{1}{n^{1/2}\kappa^{1/4}}\right)\right).\\
&\langle \hat{\bmu}_1, \bmu\rangle =\tilde \bmw_1-\tilde \bmw_2=
\frac{2c}{\beta}\sqrt{\cC(x)}
=\frac{2c }{\beta\sqrt{x-m_\phi}}\left(1+\OO_\prec\left(\frac{1}{\phi n^{1/2}\kappa^{1/4}}\right)\right).\label{e:innerp2u}
\end{align}
We also recall  from \eqref{e:mphi} that 
\begin{align}\label{e:x}
1+\left((\phi^2-1)/x-x\right)m_\phi+m_\phi^2=0, \quad \frac{m_\phi(x)}{x-m_\phi(x)}= \frac{\tau}{\phi},\quad
\tau=\frac{1}{\la^2}+\OO_\prec\left(\frac{1}{n^{1/2}\kappa^{1/4}}\right).
\end{align}
Then by taking derivative with respect to $x$ on both sides of the first expression in \eqref{e:x}, we get the following expression of $m_\phi'(x)$,
\begin{align*}
\left(-(\phi^2-1)/x^2-1\right)m_\phi+\left((\phi^2-1)/x-x\right)m_\phi'+2m_\phi m_\phi'=0,\quad
m_\phi'=\frac{((\phi^2-1)/x^2+1)m_\phi}{m_\phi-1/m_\phi}.
\end{align*}

Explicitly, we can solve for $m_\phi'(x)$ in terms of $m_\phi(x)$ as
\begin{align*}
m_\phi'&= \frac{((\phi^2-1)/x^2+1)}{1-1/m_\phi^2}
=\frac{2-(m_\phi+1/m_\phi)\frac{1}{x}}{1-1/{m_\phi}^2}
=\frac{2-(m_\phi+1/m_\phi)\frac{1}{(\phi/\tau+1)m_\phi}}{1-1/{m_\phi}^2}\\
&=\frac{2-(m_\phi+1/m_\phi)\frac{1}{(\phi/\tau+1)m_\phi}}{1-1/{m_\phi}^2}
=\left(1+\frac{\phi}{\tau}\right)^{-1}\left(1+\frac{2\phi/\tau}{1-1/m_\phi^2}\right).
\end{align*}
With the expression of $m_\phi'(x)$, we can use \eqref{e:computec} to compute $c$
\begin{align*}
\left(1+\OO_\prec\left(\frac{1}{ n^{1/2}\kappa^{3/4}}\right)\right)
&=-2c^2m_\phi \left( \frac{m_\phi'(x)}{m_\phi(x-m_\phi(x))}+\frac{(m'_\phi-1)}{(x-m_\phi)^2}\right)\\
&=-2\frac{c^2m_\phi}{(\phi/\tau)^2} \left( \frac{\phi}{\tau}\frac{m_\phi'(x)}{m_\phi^2}+\frac{(m'_\phi-1)}{m_\phi^2}\right)\\
&=-2\frac{c^2m_\phi}{(\phi/\tau)^2}\frac{2(\phi/\tau)}{m_\phi^2-1}
=4\frac{c^2m_\phi}{\phi/\tau}\frac{1}{1-m_\phi^2}.
\end{align*}
Comparing with \eqref{e:innerp}, we conclude that
\begin{align}\label{e:innerp3}
\langle \hat\bmv_1, \bmv\rangle=\frac{2c \sqrt{m_\phi}}{\la \sqrt\phi}\left(1+\OO_\prec\left(\frac{1}{n^{1/2}\kappa^{1/4}}\right)\right)
=\frac{\sqrt{(1-m_\phi^2)/\tau}}{\la}\left(1+\OO_\prec\left(\frac{1}{n^{1/2}\kappa^{3/4}}\right)\right).
\end{align}
We recall from \eqref{e:xmexp}  that
\begin{align}\label{e:mphiexx}
m_\phi=\left(1+\OO_\prec\left(\frac{1}{n^{1/2}\kappa^{1/4}}\right)\right) \sqrt{\frac{\phi+\la^2}{\la^2(1+\la^2\phi)}}.
\end{align}
and $\kappa\asymp (\la-1)^2$.
Then we conclude that
\begin{align}\label{e:v1iexx}
\langle \hat\bmv_1, \bmv_1\rangle=\left(1+\OO_\prec\left(\frac{1}{n^{1/2}(\la-1)^{3/2}}\right)\right)\sqrt{\frac{(\la^4-1)}{\la^4+\la^2/\phi}}
\sim \sqrt{1-\frac{1}{\la^4}},
\end{align}
as $\phi$ goes to infinity. 

For the inner product $\langle \hat\bmu_1, \bmu\rangle$ from \eqref{e:innerp2u}, we have
\begin{align*}
&\phantom{{}={}}\langle \hat{\bmu}_1, \bmu\rangle 
=\langle \hat{\bmv}_1, \bmv\rangle\sqrt{\frac{1}{m_\phi}\frac{m_\phi}{x-m_\phi}}\left(1+\OO_\prec\left(\frac{1}{n^{1/2}\kappa^{3/4}}\right)\right)\\
&=\sqrt{\frac{\la^4-1}{\la^2(\la^2+\phi)}}\left(1+\OO_\prec\left(\frac{1}{n^{1/2}\kappa^{3/4}}\right)\right)
=\left(1+\OO_\prec\left(\frac{1}{n^{1/2}(\la-1)^{3/2}}\right)\right)\sqrt{\frac{\la^4-1}{\la^2(\la^2+\phi)}},
\end{align*}
where for the second line, we plugged in the estimates \eqref{e:x}, \eqref{e:innerp3} and \eqref{e:mphiexx}.
This finishes the proof of \eqref{e:strong1} and \eqref{e:strong22}.
\end{proof}

In the following, we study the case when $\la\leq 1+n^{-2/3+\fc}$ and prove \eqref{e:weak1}. 
\begin{proof}[Proof of \eqref{e:weak1}]
We recall the Hermitized matrix from \eqref{e:hermitc}. Its nonzero eigenvalues are given by $\hat s_1, -\hat s_1, \hat s_2, -\hat s_2,\cdots, \hat s_n, -\hat s_n$, and the corresponding normalized eigenvectors are given by 
$[\hat\bmv_1, \hat\bmu_1]/\sqrt 2, [\hat\bmv_1, -\hat\bmu_1]/\sqrt 2, [\hat\bmv_2, \hat\bmu_2]/\sqrt 2, [\hat\bmv_2, -\hat\bmu_2]/\sqrt 2,\cdots, [\hat\bmv_n, \hat\bmu_n]/\sqrt 2, [\hat\bmv_n, -\hat\bmu_n]/\sqrt 2$. Then its Green's function is given by
\begin{align}\begin{split}\label{e:hermitcc}
&\phantom{{}={}}\left(z-\left[
\begin{array}{cc}
\bmv/\sqrt 2 & \bmv\sqrt{2}\\
\bmu/\sqrt 2 & -\bmu\sqrt{2}
\end{array}
\right]
\left[
\begin{array}{cc}
\beta & 0\\
0  & -\beta
\end{array}
\right]
\left[
\begin{array}{cc}
\bmv/\sqrt 2 & \bmv\sqrt{2}\\
\bmu/\sqrt 2 & -\bmu\sqrt{2}
\end{array}
\right]^\top
+\left[
\begin{array}{cc}
0 & Z\\
Z^\top & 0
\end{array}
\right]\right)^2\\
&=\sum_{i=1}^n \frac{[\hat \bmv_i, \hat \bmu_i][\hat \bmv_i, \hat \bmu_i]^\top}{2(z-\hat s_i)} + \frac{[\hat \bmv_i, -\hat \bmu_i][\hat \bmv_i, -\hat \bmu_i]^\top}{2(z-\hat s_i)}+\frac{WW^\top}{z}\\
&=\left(I+A(z)(I-U^\top A(z))^{-1}U^\top\right)G(z),
\end{split}\end{align}
where 
\begin{align*}
A(z)=G(z)\left[
\begin{array}{cc}
\bmv/\sqrt 2 & \bmv\sqrt{2}\\
\bmu/\sqrt 2 & -\bmu\sqrt{2}
\end{array}
\right]
\left[
\begin{array}{cc}
\beta & 0\\
0  & -\beta
\end{array}
\right],\quad
U=
\left[
\begin{array}{cc}
\bmv/\sqrt 2 & \bmv\sqrt{2}\\
\bmu/\sqrt 2 & -\bmu\sqrt{2}
\end{array}
\right],
\end{align*}
are defined as in \eqref{e:AU}, and $W$ is an $(m+n)\times(m-n)$ matrix with columns the eigenvectors of the Hermitized matrix \eqref{e:hermitc} corresponding to eigenvalue $0$.

We conjugate \eqref{e:hermitcc} by $U$ on both sides, and get a $2\times 2$ matrix
\begin{align*}
U^\top\left(I+A(z)(I-U^\top A(z))^{-1}U^\top\right)G(z)U
=(I-U^\top A(z))^{-1}(U^\top G(z)U).
\end{align*}
This matrix contains the information of projection of $[\bmv, \pm \bmu]$ in the directions of $[\hat \bmv_i, \pm \hat \bmu_i]$. More precisely
\begin{align}\begin{split}\label{e:decomp}
&((I-U^\top A(z))^{-1}(U^\top G(z)U))_{11}\\
&=\sum_{i=1}^n \frac{\langle [\hat \bmv_i, \hat \bmu_i], [\bmv, \bmu]\rangle^2}{4(z-\hat s_i)} + \frac{\langle [\hat \bmv_i, -\hat \bmu_i], [\bmv, \bmu]\rangle^2}{4(z-\hat s_i)}+\frac{\|W[\bmv, \bmu]\|_2^2}{2z},\\
&((I-U^\top A(z))^{-1}(U^\top G(z)U))_{22}\\
&=\sum_{i=1}^n  \frac{\langle [\hat \bmv_i, \hat \bmu_i], [\bmv, -\bmu]\rangle^2}{4(z-\hat s_i)} + \frac{\langle [\hat \bmv_i, -\hat \bmu_i], [\bmv, -\bmu]\rangle^2}{4(z-\hat s_i)}+\frac{\|W[\bmv, -\bmu]\|_2^2}{2z}.
\end{split}\end{align}
We notice that since $z\in \bC^+$, the imaginary part of each term on the righthand side of \eqref{e:decomp} is negative. By taking the sum of the two terms in \eqref{e:decomp} and imaginary part on both sides, we get
\begin{align}\label{e:evbb}
|\Im[\Tr ((I-U^\top A(z))^{-1}(U^\top G(z)U))]|
\geq \frac{\Im[z] (\langle \hat \bmv_1, \bmv\rangle^2+\langle \hat \bmu_1, \bmu\rangle^2)}{|z-\hat s_1|^2}.
\end{align}
By taking $z=\hat s_1+\ri \eta$ in \eqref{e:evbb}, we get the upper bound
\begin{align}\label{e:evub}
 \langle \hat \bmv_1, \bmv\rangle^2+\langle \hat \bmu_1, \bmu\rangle^2\leq \eta |\Im[\Tr ((I-U^\top A(z))^{-1}(U^\top G(z)U))]|.
\end{align}

In the following, we estimate the lefthand side of \eqref{e:evbb}. We recall $\cA(z), \cB(z), \cC(z)$ from \eqref{e:defABC}, they are well defined for $z\in \bC^+$, and Theorem \ref{t:isotropiclaw} implies
\begin{align}\begin{split}\label{e:ABCest}
&\cA(z)=m_\phi(z)+\OO_\prec\left(\sqrt{\frac{|\Im[m_\phi(z)]|}{n\eta}}+\frac{1}{n\eta}\right),\\
&\cB(z)=\OO_\prec\left(\frac{1}{\phi}\left(\sqrt{\frac{|\Im[m_\phi(z)]|}{n\eta}}+\frac{1}{n\eta}\right)\right),\\
&\cC(z)=\frac{1}{z-m_\phi(z)}+\OO_\prec\left(\frac{1}{\phi^2}\left(\sqrt{\frac{|\Im[m_\phi(z)]|}{n\eta}}+\frac{1}{n\eta}\right)\right),
\end{split}\end{align}
The matrix $U^\top G(z) U$ can be expressed in terms of $\cA(z), \cB(z), \cC(z)$
\begin{align}\label{e:UGU}
U^\top G(z)U=\frac{1}{2}\left[
\begin{array}{cc}
\cA(z)+2\cB(z)+\cC(z) & \cA(z)-\cC(z)\\
\cA(z)-\cC(z) & \cA(z)-2\cB(z)+\cC(z)
\end{array}
\right]. 
\end{align}
$(I-U^\top A(z))$ is a $2\times 2$ matrix, we can invert it by Cramer's rule 
\begin{align}\label{e:Ainv}
\frac{1}{\det(I-U^\top A(z))}
\left[
\begin{array}{cc}
1+\beta(\cA(z)-2\cB(z)+\cC(z))/2 & \beta(\cA(z)-\cC(z))/2\\
-\beta(\cA(z)-\cC(z))/2 & 1-\beta(\cA(z)+2\cB(z)+\cC(z))/2
\end{array}
\right],
\end{align}
and the determinant is given by
\begin{align}\label{e:deta}
\det(I-U^\top A(z))=\beta^2\left(\left(\frac{1}{\beta}-\cB(z)\right)^2-\cA(z)\cC(z)\right)
\end{align}
By plugging \eqref{e:UGU} and \eqref{e:Ainv} into \eqref{e:evbb}, we get
\begin{align}\label{e:tra}
\Tr ((I-U^\top A(z))^{-1}(U^\top G(z)U))
=\frac{\cA(z)+\cC(z)}{\det(I-U^\top A(z))}.
\end{align}
To use \eqref{e:evub}, we will take $z\in \bC^+$ in a small neighborhood of the spectral edge $\phi+1$. Let $z=\phi+1+\kappa+\ri\eta$, where $\kappa, \eta\ll1$. Then thanks to the explicit formula of $m_\phi(z)$ from \eqref{e:mpsi}, in this region, we have
\begin{align*}
m_\phi(z)=1-C\sqrt{\kappa+\ri\eta}+\OO(|\kappa+\ri\eta|),
\end{align*}
for some constant $C>0$, and 
\begin{align}\label{e:AC1bb}
\cA(z)+\cC(z)=1-(C+C/\phi)\sqrt{\kappa+\ri\eta}+\OO_\prec\left(\sqrt{\frac{|\Im[m_\phi(z)]|}{n\eta}}+\frac{1}{n\eta}+|\kappa+\ri\eta|\right)\asymp 1.
\end{align}
Recall that we have $\beta=\la\sqrt \phi$. 
For the denominator in \eqref{e:tra}, using \eqref{e:deta} and the estimates \eqref{e:ABCest}, we get
\begin{align}\begin{split}\label{e:dttIUA}
&\phantom{{}={}}\det(I-U^\top A(z))
=1-\frac{\beta^2 m_\phi(z)}{z-m_\phi(z)}+\OO_\prec\left(\sqrt{\frac{|\Im[m_\phi(z)]|}{n\eta}}+\frac{1}{n\eta}\right)\\
&=(1-\la^2)+\la^2(C+C/\phi)\sqrt{\kappa+\ri\eta}+\OO_\prec\left(\sqrt{\frac{|\Im[m_\phi(z)]|}{n\eta}}+\frac{1}{n\eta}+|\kappa+\ri\eta|\right).
\end{split}\end{align}
There are two cases, either $\la\in[1-n^{-1/3+2\fc}, 1+n^{-1/3+\fc}]$, or $\la\leq 1-n^{-1/3+2\fc}$. If $\la\in[1-n^{-1/3+2\fc}, 1+n^{-1/3+\fc}]$, we can take $\kappa\in [-n^{-2/3+3\fc}, n^{-2/3+3\fc}]$, and $\eta=n^{-2/3+6\fc}$, then 
\begin{align}\begin{split}\label{e:ddertt}
&\phantom{{}={}}|\det(I-U^\top A(z))|\geq |\Im[\det(I-U^\top A(z))]|\\
&=\la^2(C+C/\phi)\Im[\sqrt{\kappa+\ri\eta}]+\OO_\prec\left(\sqrt{\frac{|\Im[m_\phi(z)]|}{n\eta}}+\frac{1}{n\eta}+|\kappa+\ri\eta|+|1-\la|\right)\gtrsim \sqrt\eta .
\end{split}\end{align}
Then \eqref{e:evub},\eqref{e:tra}, \eqref{e:AC1bb} and \eqref{e:ddertt}  imply that
\begin{align}\label{e:v1v}
 \langle \hat \bmv_1, \bmv\rangle^2+\langle \hat \bmu_1, \bmu\rangle^2\lesssim \sqrt\eta=n^{-1/3+3\fc}.
\end{align}

If $\la\leq 1-n^{-1/3+2\fc}$, we take $\kappa\in [-n^{-2/3+3\fc}, n^{-2/3+3\fc}]$ and $\eta=n^{\fc-1}/(1-\la)$ in \eqref{e:dttIUA}
\begin{align}\begin{split}\label{e:dtIIH}
|\det(I-U^\top A(z))|
&\geq 1-\la^2+\OO_\prec\left(\sqrt{\frac{|\Im[m_\phi(z)]|}{n\eta}}+\frac{1}{n\eta}+\sqrt{|\kappa+\ri\eta|}\right)\\
&\geq 1-\la^2+\OO_\prec\left(\sqrt{\frac{\sqrt{|\kappa|+\eta}}{n\eta}}+\frac{1}{n\eta}+\sqrt{|\kappa|+\eta}\right)\\
&\gtrsim 1-\la^2+\OO_\prec\left(\frac{1}{n\eta}+\sqrt{|\kappa|+\eta}\right)
\gtrsim 1-\la^2.
\end{split}\end{align}
Then \eqref{e:evub},\eqref{e:tra}, \eqref{e:AC1bb} and \eqref{e:dtIIH} imply that
\begin{align}\label{e:v2v}
 \langle \hat \bmv_1, \bmv\rangle^2+\langle \hat \bmu_1, \bmu\rangle^2\lesssim \eta|\Tr ((I-U^\top A(z))^{-1}(U^\top G(z)U))|\lesssim \frac{n^{\fc}}{n(1-\la)^2}.
\end{align}
Then \eqref{e:v1v} and \eqref{e:v2v} together imply that
\begin{align*}
 \langle \hat \bmv_1, \bmv\rangle^2+\langle \hat \bmu_1, \bmu\rangle^2\lesssim n^{8\fc} \min\left\{n^{-1/3}, \frac{1}{n(1-\la)^2}\right\}.
\end{align*}
This finishes the proof of \eqref{e:weak1}.
\end{proof}

Theorem \eqref{t:ev} gives the behavior of the projection of the singular vector $\hat\bmv_1$ associated with the largest singular value of $\beta \bmv \bmu^\top + Z$ on the signal direction. At the critical value $\lambda=1$, it states
\begin{align*}
|\langle \hat\bmv_1, \bmv\rangle|^2\lesssim n^{-1/3+\oo(1)}.
\end{align*}
We believe that it is optimal up to the $\oo(1)$ error in the exponent. 
More precisely, we conjecture that
exactly at the critical signal strength, $\beta=\sqrt{\phi}$, the projection of the singular vector $\hat\bmv$ associated with the largest singular value of $\beta \bmv\otimes \bmu + Z$ on the signal $\bmv$ direction satisfies
\begin{align}\label{e:converge}
n^{1/3}|\langle \hat\bmv, \bmv\rangle|^2\rightarrow \Theta,
\end{align}
where $\Theta$ is a random variable of size $\OO(1)$.
The above statement \eqref{e:converge} for low rank perturbations of Gaussian unitary matrices have been proven in \cite{bao2020eigenvector}, where they give explicit characterization of the limiting objection $\Theta$.

\small
\bibliographystyle{abbrv}
\bibliography{References}

\end{document}